\documentclass[11pt]{article}

\usepackage[latin1]{inputenc}
\usepackage{color}

\usepackage{amsmath}
\usepackage{amssymb}
\usepackage{graphicx}

\newcommand{\eps}{\epsilon}
\newcommand{\qed}{$\square$}

\newcommand{\al}{\alpha}

\newcommand{\ep}{l}

\def\lr#1{\langle{#1}\rangle}

\newcommand{\NN}{\mathbb{N}}
\newcommand{\RR}{\mathbb{R}}
\newcommand{\CC}{\mathbb{C}}

\newcommand{\ut}{\underline t}

\newtheorem{theorem}{Theorem}[section]

\newtheorem{remark}{Remark}[section]

\newtheorem{hypothesis}{Hypothesis}[section]

\newtheorem{definition}{Definition}[section]

\newtheorem{example}{Example}[section]

\newtheorem{proposition}[theorem]{Proposition}

\newtheorem{corollary}[theorem]{Corollary}

\newtheorem{lemma}[theorem]{Lemma}

\numberwithin{equation}{section}

\begin{document}

\begin{center}
\Large
\bf
Weakly Hyperbolic Systems
 by Symmetrization
\end{center}

\begin{center}
F. Colombini 
\footnote{Dipartimento di Matematica, Universit\`a di Pisa, Pisa, Italy.
{\tt  ferruccio.colombini@unipi.it}
}
\quad
T. Nishitani
\footnote{Department of Mathematics,
Graduate School of Science, 
Osaka University, Osaka, Japan. 
{\tt nishitani@math.sci.osaka-u.ac.jp}
}
\quad
 and 
 \quad
 J.Rauch \footnote{Department of Mathematics, University of Michigan, Ann Arbor, Michigan, USA.
 {\tt rauch@umich.edu}}
\end{center}

\begin{abstract}  We prove Gevrey well posedness of the Cauchy
problem for general linear systems whose principal symbol is
hyperbolic and  coefficients are 
sufficiently Gevrey regular in $x$ and either lipschitzian or h\"olderian
in time.  
Such
results date to the seminal paper of Bronshtein.
The proof is by an energy method using a pseudodifferential
symmetrizer.  The  construction of the symmetrizer 
is based on  a lyapunov function for
ordinary differential equations.   The method yields 
new estimates and 
existence uniformly for spectral truncations and parabolic
regularizations.  
\end{abstract}

{\bf MSC Classification.}  35L45, 35L40

{\bf Keywords.}  Weakly hyperbolic,  Gevrey regularity, symmetrizer,  energy method,
characteristics of variable multiplicity.

\section{Introduction}

A partial differential operator in $t,x\in \RR^{1+d}$
is called {\bf hyperbolic} when $t=0$ is non characteristic  and 
the characteristic polynomial has only real roots $\tau$ for 
arbitrary $\xi\in \RR^d\setminus 0$. For the first order systems
that we consider,
\begin{equation}
\label{eq:Cauchy}
Lu \ =\ \partial_tu-\sum_{j=1}^dA_j(t,x)\partial_{x_j}u+B(t,x)u\ =\  f,
\quad
u(0,\cdot) = g,
\end{equation}
the coefficients $A_j$ and $B$ are $m\times m$ matrix valued.
The principal symbol is 
$$
A(t,x,\xi) \ :=\ \sum_{j=1}^d A_j(t,x)\,\xi_j\ .
$$
Hyperbolicity, assumed throughout, means
\begin{equation}
\label{eq:hyperbolicity}
\forall t,x,\xi\in \RR\times \RR^d\times\RR^d,
\qquad
{\rm Spectrum} \, A(t,x,\xi)\ \subset\ \RR\,.
\end{equation}
For the noncharacteristic Cauchy problem,
\eqref{eq:hyperbolicity}
is
a necessary condition for the Cauchy problem
to be well set for non analytic data.
The condition is not sufficient for well posedness for 
$C^\infty$ data.  For most non strictly hyperbolic
scalar operators, most lower order terms lead to initial value
problems that are ill posed in the $C^\infty$ category above.
The generic ill posedness holds even if the coefficients are 
real analytic functions or even constant.

For real analytic hyperbolic operators,
\cite{Ivr} and \cite{Trep} showed that for Gevrey initial data,
${\cal G}^s$ with 
$1<s<s_0$,
there are Gevrey solutions.
No conditions of E.Levi type 
 is needed.  It came as a surprise to 
many, including us, when
Bronshtein \cite{Bronsh:2} proved that the Cauchy problem for linear
 hyperbolic partial differential operators whose coefficients are finitely smooth in time 
 and Gevrey in $x$ is well posed for Gevrey data. 
 Bronshtein, Ohya-Tarama \cite{OhTa}, and Kajitani \cite{Kaji}, \cite{Kaji:2}
  used  parametrix constructions  either by examinig the resolvent close the imaginary
  axis or by Fourier Integral Operator 
  constructions.
     The papers
     \cite{CDS}, \cite{CJS}, \cite{Ni:1}, \cite{Ni:2} 
     use energy methods of increasing complexity.
   In this paper we 
   introduce an  energy method that we think is as simple
   as the very simplest of these and also very natural.
  Our estimates are proved while ignoring the detailed behavior of the eigenvalue
  crossings. We call this as working with our eyes shut.

   A standard approach to  proving well posedness for Gevrey data for  hyperbolic systems is to multiply  the reference system by the operator of cofactor symbol
   to  reduce  to  scalar operators.
   That approach has at least two defects.  First 
   applying the cofactor matrix requires 
   that the coefficients have a number of derivatives in time roughly equal to the size of the matrix.
   Second, this totally ignores the  system structure.   For example if a system is merely
   two copies of a strictly hyperbolic system, the cofactor  approach  immediately replaces
   the problem with one that is much less well behaved.
   
   We study first order  hyperbolic systems and 
   prove Gevery well posedness,  by proving {\it a priori} estimates by constructing a 
    pseudodifferential symmetrizer.
   The symmetrizer is  motivated by a special Lyapunov function 
   for asymptotically stable 
   constant coefficient first order systems of  linear 
   ordinary differential equations. The proof not only gives straightforward  {\it a priori} estimates,
 but also clarifies some effects coming from the block structure of the system.
 It does not at all 
 look closely at the eigenvalue crossings and that is its principal strength.
 
This paper  discusses only the existence and uniqueness of solutions.  
The method of  
 \cite{CoRau} 
gives the natural precise estimate for the influence domain.
In particular this allows one to eliminate our hypothesis that the coefficients are 
independent of $x$ outside a compact subset of space.

To our systems  we
associate, in  Hypothesis \ref{hyp:theta},  an  index $0\le \theta\le m-1$.
The value of 
 $\theta$ measures roughly whether the   Taylor polynomial of degree 
$N=\max\{2\theta, m\}$
of the symbol can be uniformly block diagonalized with blocks of size $\theta+1$.
It is always satisfied with  $\theta=m-1$.   

The functions uniformly Gevrey $s$ on $\RR^d$ are denoted ${\mathcal G}^{s}(\RR^d)$
and those of compact  support by ${\mathcal G}^{s}_0(\RR^d)$.
In  the results below, the Gevrey index $s_0$
is cruder, that is smaller, than the sharp results of 
\cite{CDS}, \cite{CJS} valid in special cases.
The result for coefficients
lipshitzian in time is the following.

\begin{theorem}
\label{thm:bronshtein1}
Suppose Hypothesis \ref{hyp:theta} is satisfied.
Define
$$
s_0\ :=\
\max
\Big\{
    \frac{ 2+6\theta }
    {1+6\theta}
    \,,\,
    \frac{ 3+4\theta }
    { 2+4\theta }
    \Big\}\,.
    $$
For some $1<s\le s_0$ 
suppose the coefficients 
$A_j(t,x)$ (resp. $B(t,x)$) are lipschitzian
(resp. continuous) in time uniformly on compact sets with 
values in the elements of ${\mathcal G}^{s}(\RR^d)$ 
that are constant outside a fixed compact set, 
$g\in {\mathcal G}_0^{s}(\RR^d)$, and $f\in L^1_{loc}(\RR\,; {\mathcal G}_0^{s}(\RR^d))$.
  Then there is a $T_0>0$ and a unique 
local solution
$
u\in C([0,T_0]\,; {\mathcal G}_0^{s}(\RR^d))
$
to the Cauchy problem 
\eqref{eq:Cauchy}.
\end{theorem}

\begin{remark}
\label{rem:Thm1.1}
\rm
The proof shows in addition that 
for all constants $c>0$ and $T>0$ the interval of existence can be chosen uniformly for the data satisfying
\[
\int |{\hat g}(\xi)|^2e^{c\lr{\xi}^{1/s}}d\xi
\ +\ 
\int_0^T\Big(\int|{\hat f}(t,\xi)|^2e^{c\lr{\xi}^{1/s}}d\xi\Big)^{1/2}
\ dt\ <\ \infty\,.
\]
An analogous remark applies
Theorem \ref{thm:lowreg}.
\end{remark}

The next result concerns
equations with coefficients H\"older continuous of order $\kappa$
in  time.  

\begin{theorem}
\label{thm:lowreg}
Suppose that $0<\kappa<1$ and that Hypothesis \ref{hyp:theta} holds.  Define
$$
 s_0\ :=\
  \min{
  \bigg\{
  \frac{2+3\theta}{2+3\theta-\kappa}
  \ ,\
  \max{
  \Big\{
  \frac{2+6\theta}{1+6\theta}
  \, ,\,
  \frac{3+4\theta}{2+4\theta}
  \Big\}
  }
  \bigg\}
  }
    \,.
    $$
Suppose that the map 
the  $t\mapsto A_j(t,\cdot)$ 
(resp. $t\mapsto B(t,\cdot)$) is $\kappa$ H\"older continuous 
(resp. continuous)
  in time uniformly on compact sets with 
values in the elements of ${\mathcal G}^{s}(\RR^d)$ that are constant outside a fixed compact 
subset of $\RR^d$,
$g\in {\mathcal G}_0^{s}(\RR^d)$, and $f\in L^1_{loc}(\RR\,; {\mathcal G}_0^{s}(\RR^d))$.
     Then the conclusion of Theorem \ref{thm:bronshtein1} holds. 
\end{theorem}

The idea of the symmetrization is easy. We multiply by a positive hermitian
pseudodifferential operator to derive estimates.
  The change of variables $v=e^{a\langle D \rangle^{\rho} t}u$
replaces the operator $L$ by $L-a\langle D\rangle^{\rho}$.   Choosing $a>>1$ and $0<\rho<1$
appropriately,
  the matrix
$$
M(t,x,\xi) \ =\ 
 A(t,x,\xi) \ +\ 
  B(t,x) \ -\
  a\langle\xi\rangle^\rho
$$
 has spectrum with real part 
$\le -\langle\xi\rangle^{\rho}$
for all $t,x,\xi$.  For the ordinary differential equation $X^\prime=MX$, the positive definite matrix 
$$
R(t,x,\xi)\ :=\
\int_0^\infty
\big(e^{Ms}\big)^*
e^{Ms}\ ds
$$
defines a strict lyapunov function, that is $RM +M^* R <0$.
Our symmetrizer is based on  
 $R(t,x,D)$.  
 This multiplier method has many advantages.  For example, it yields
estimates uniform in $\eps$ for the regularized operators
$$
\partial_t + \chi(\eps D)\Big(\sum_jA_j\partial_j + B
\Big)\chi(\eps D)\,,
\qquad
\chi\in {\cal S}(\RR^d),\quad
\chi(0)=1\,,
$$
as well as for parabolic regularizations,
$$
\partial_t \ +\  
\sum_jA_j\partial_j 
\ +\ B
 \ -\ 
 \eps\,\Delta\,.
$$
The first is used to prove existence and is related to the spectral
method analysed in \cite{CR2}.  
This yields two more ways that these   very weakly hyperbolic
equations
are
in line with other
hyperbolic Cauchy problems.

\section{Gevrey operators}

\subsection{Symbol classes and conjugation}

Denote
\begin{equation}
\label{eq:defanglexi}
\lr{\xi}_{\ell}\ :=\
\sqrt{\ell^2+|\xi|^2}
\ =\ 
\ell\sqrt{1\, +\, |\xi/\ell|^2}
\end{equation}
where $\ell\geq 1$ is a positive large parameter. We write $\lr{\xi}_1=\lr{\xi}$ and note that $\lr{\xi}\leq \lr{\xi}_{\ell}\leq \ell\lr{\xi}$.
\begin{definition}
\label{dfn:gs}
If $1<s<\infty$, the function $a(x)\in  C^{\infty}(\RR^d)$ belongs  to  ${\mathcal G}^{s}(\RR^d)$ if  there exist $C>0$, $A>0$ such that
\[
\forall  x\in {\mathbb R}^d,\quad
\forall  \alpha\in\NN^d,
\qquad
|\partial_x^{\al}a(x)|
\ \leq\ C A^{|\al|}|\al|!^s\,.
\]
Denote ${\mathcal G}_0^{s}(\RR^d):={\mathcal G}^{s}(\RR^d)\cap C_0^{\infty}(\RR^d)$.
\end{definition}
\begin{definition}
\label{dfn:rhodelta} 
For $0<\delta\le \rho\le 1$, the function  $a(x,\xi;\ell)\in C^{\infty}(\RR^d\times\RR^d)$ 
belongs to  ${\tilde S}^m_{\rho,\delta}$ if for all $\al$, $\beta\in \NN^d$ there is $C_{\al\beta}$ independent of $\ell,x,\xi$ such that
\[
\big|
\partial_x^{\beta}\partial_{\xi}^{\al}a(x,\xi;\ell)
\big|
\ \leq\
 C_{\al\beta}\lr{\xi}_\ell^{m-\rho|\al|+\delta|\beta|}.
\]
Denote  ${\tilde S}^m:={\tilde S}^m_{1,0}$.
\end{definition}

\begin{definition}
\label{dfn:kurasu}
For  $1<s$, $m\in\RR$, 
the function $a(x,\xi;\ell)\in C^{\infty}(\RR^d\times\RR^d)$ 
belongs to  ${\tilde S}_{(s)}^m$ if there exist $C>0$, $A>0$ independent of $\ell,x,\xi$ such that for all $\alpha$, $\beta\in\NN^d$,
\[
\big|
\partial_x^{\beta}\partial_{\xi}^{\alpha}a(x,\xi;\ell)
\big|
\ \leq\ 
C\,A^{|\alpha+\beta|}\ 
|\alpha+\beta|!^s\
\lr{\xi}_{\ell}^{m-|\alpha|}
\,.
\]
\end{definition}

We often write $a(x,\xi)$ for $a(x,\xi,\ell)$ dropping the $\ell$. If $a(x,\xi)$ is the symbol of a differential operator of order $m$ with coefficients $a_{\alpha}(x)\in {\mathcal G}^{s}(\RR^d)$ 
then $a(x,\xi)\in {\tilde S}^m_{(s)}$ because $|\partial_{\xi}^{\beta}\xi^{\alpha}|\leq CA^{|\beta|}|\beta|!\lr{\xi}_{\ell}^{|\alpha|-|\beta|}$ and $|\partial_x^{\beta}a_{\alpha}(x)|\leq C_{\alpha}A_{\alpha}^{|\beta|}|\beta|!^s$ for any $\beta\in\NN^d$. 
\begin{proposition}
\label{pro:weyl:1} 
Suppose $1/2\leq \rho<1$, $s=1/\rho$, and
$a(x,\xi)$ be $m\times m$ matrix valued with entries in $ {\tilde S}_{(s)}^1$ 
and $\partial_x^{\alpha}a(x,\xi)=0$ outside $|x|\leq R$ with some $R>0$ if $|\alpha|>0$. Then  
the operator $b(x,D)=e^{\tau\lr{D}_{\ell}^{\rho}}a(x,D)e^{-\tau\lr{D}_{\ell}^{\rho}}$ is
for small $|\tau|$
a pseudodifferential operator with symbol  given by
\[
b(x,\xi)\ =\ 
\sum_{|\alpha|\leq  k}\frac{1}{\alpha !}\ D_x^{\alpha}a(x,\xi)\
(\tau\nabla_{\xi}\lr{\xi}_{\ell}^{\rho})^{\alpha}
\ +\ R(x,\xi)
\]
with 
$R(x,\xi)\in {\tilde S}^{\,\max{\{\rho-k(1-\rho)\,,\,-1+\rho\}}}$  and $D_{x_j}:=-i\partial/\partial x_j$. 
\end{proposition}
Proposition \ref{pro:weyl:1} is not new.  For completeness
a  proof is given in \S \ref{Sec:composition}.

If $a$ were real analytic 
in $x$ then the sum on  the right would be 
$$
\sum_{|\alpha|\leq  k}\frac{\partial_x^{\alpha}a}{\alpha !} \ (-iy)^\alpha
\ = \ 
a(x-iy\,,\,\xi)\ +\ O(|y|^{k+1})\,,
\qquad
y=\tau\nabla_\xi\langle\xi\rangle^\rho_\ell
\,.
$$
For large $\xi$, $y$ tends to zero because $\rho<1$.
Therefore, this is a very small displacement in
the complex direction.   In the early work of \cite{Ivr}, \cite{Trep}
the coefficients were analytic and one could make such complex
displacements.  For our problems, the coefficients are not analytic and
the replacement for complex displacement is to put complex arguments into
Taylor polynomials.  An alternative strategy is to take
an almost analytic extension of $a$ that satisfies the Cauchy-Riemann
equations with error $O(|y|^\infty)$ at $y=0$.
\begin{corollary}
\label{cor:esurei} If $a(x,\xi)\in {\tilde S}_{(s)}^0$ then $e^{\tau\lr{D}_{\ell}^{\rho}}a(x,D)e^{-\tau\lr{D}_{\ell}^{\rho}}\in {\rm Op}\,{\tilde S}^{0}$ for small $|\tau|$.
\end{corollary}
{\bf Proof:} Choose $k$ so that $\rho-k(1-\rho)\leq 0$. 
Then $D_x^{\alpha}a(x,\xi)(\tau\nabla_{\xi}\lr{\xi}_{\ell}^{\rho})^{\alpha}\in {\tilde S}^{-(1-\rho)|\al|}\subset {\tilde S}^0$ for $a\in {\tilde S}^0$. The assertion follows from Proposition \ref{pro:weyl:1}.
\hfill
\qed

\subsection{The block size  barometer $\mathbf \theta$}

Introduce an integer valued  parameter $0\le \theta \le m-1$ that measures the extent to which
the principal symbol can be block diagonalized by matrices bounded with bounded 
inverse.   For example in the strictly hyperbolic case, blocks of size 1 are attainable.
By convention $\theta$ is one smaller than the block size.   Block size $m$ and 
therefore $\theta=m-1$ is always possible. The definition of 
$\theta$ is not directly given in these terms.   The relation to block size
is discussed in the examples below.

Assume that $A_j(t,x)\in C^0(\RR\,;C^{\infty}(\RR^d))$ and all eigenvalues of $A(t,x,\xi)=\sum_{j=1}^dA_j(t,x)\xi_j$ are real for any $(t,x)\in\RR\times \RR^d$ and $\xi\in\RR^d$. 
From Proposition \ref{pro:hypspect} for any $T>0$ and compact set $K\subset \RR^d$ there exist $\delta>0$ and $c>0$ such that if $\zeta$ is an eigenvalue of
\begin{equation}
\label{eq:saji}
\sum_{|\al+\beta|\leq m}\frac{(is)^{|\al+\beta|}}{\al!\beta!}\partial_x^{\al}\partial_{\xi}^{\beta}A(t,x,\xi)y^{\al}\eta^{\beta}
\end{equation}
then $|{\mathsf{Im}}\,\zeta|\leq c\,|s|$ for any $|(y,\eta)|\leq 1$, $x\in K$, $|\xi|\leq 1$, $|t|\leq T$.
Define
\[
{\mathcal H}_r(t,x,\xi;\epsilon)
\ :=\
\sum_{|\al|\leq r}\frac{\epsilon^{|\al|}}{\al!}D_x^{\al}A(t,x,\xi)\xi^{\al}.
\]
Choosing $(y,\eta)=(\xi,0)$ in \eqref{eq:saji} we see that there is $\epsilon_0>0$, $c>0$ such that
\begin{equation}
\label{eq:speccalH}
\zeta\mbox{~is an eigenvalue of~}{\mathcal H}_m(t,x,\xi;\epsilon)
\ \ \Longrightarrow\ \
 |{\mathsf{Im}}\,\zeta|\leq c\,|\epsilon|
\end{equation}
for any $x\in K$, $|\xi|\leq 1$, $|\epsilon|\leq \epsilon_0$, $|t|\leq T$.
Introduce the following hypothesis.

\begin{hypothesis} 
\label{hyp:theta}
Assume the system is  $\theta${\bf -regular} with integer 
$0\le \theta\le m-1$
in the sense that 
 for any $T>0$ and any compact $K\subset \RR^d$ 
there exist  $C>0$, $c>0$ and $\epsilon_0>0$ such that with $N=\max\{2\theta,m\}$
\begin{equation}
\label{eq:Matexp}
\frac{\eps^\theta}
{C\, e^{cs\epsilon}}
\ \leq\  
 \big\| e^{is{\mathcal H_N(t,x,\xi;\epsilon)}}   \big \|
 \ \leq\
\frac{  C\, e^{cs\epsilon}}
{\epsilon^{\theta}} \,,
  \end{equation}
 for all $s\geq 0$, $0<\epsilon\leq \epsilon_0$, $|\xi| = 1$, $x\in K$, $|t|\leq T$.
 \end{hypothesis}

A system that is  $\theta$-regular is $\phi$-regular for all $\theta<\phi\le m-1$.

Denote
\[
H_N(\rho,\ell,\tau, t,x,\xi)\ :=\
\sum_{|\alpha|\leq  N}\frac{1}{\alpha !}\,
D_x^{\alpha}A(t,x,\xi)\,
\big(\tau \nabla_{\xi}\lr{\xi}_{\ell}^{\rho}\big)^{\alpha} \,.
\]
The definition of ${\mathcal H}_N$ implies that 
\[
H_N(\rho,\ell,\tau,t,x,\xi)
\ =\ 
\lr{\xi}_\ell\ {\mathcal H}_N\Big( t\,,\, x\,,\, \xi/\lr{\xi}_\ell\,;\, \tau\rho\lr{\xi}_\ell^{\rho-1} \Big)\,.
\]
Choosing $s\lr{\xi}_\ell$, $\tau\rho\lr{\xi}_\ell^{\rho-1}$ ($\tau>0$), $\xi/\lr{\xi}_\ell$ for $s$, $\epsilon$, $\xi$ in \eqref{eq:Matexp} 
yields
\begin{equation}
\label{eq:expH}
\begin{split}
\frac{ \tau^{\theta}  }
{C\, \lr{\xi}_\ell^{\theta(1-\rho)} \, e^{cs\tau\lr{\xi}_\ell^{\rho}}}
\ \leq\ 
\big\| e^{isH_N(\ell,\tau,t,x,\xi)} \big\|
\ \leq\
\frac{C\, \lr{\xi}_\ell^{\theta(1-\rho)}   \, e^{cs\tau\lr{\xi}_\ell^{\rho}}  }
{\tau^{\theta}   }
\end{split}
\end{equation}
for $|t|\leq T$, $\ell\geq \ell_0$ where $\tau, \ell_0$ are constrained to 
satisfy 
\begin{equation}
\label{eq:constraint6}
0\ <\ \tau \ell_0^{\rho-1}
\ \leq\  \epsilon_0\,.
\end{equation}

{\begin{example}
\label{example:1}
Estimate  \eqref{eq:Matexp} always holds with $\theta=m-1$.
\rm  
Indeed write ${\mathcal H}_N={\mathcal H}_m+L_N$ where $\|L_N\|\leq C\epsilon^{m+1}$. Take an orthogonal matrix $T$ such that $T{\mathcal H}_mT^{-1}$ to be upper triangular. Let $S={\rm diag}(1,\epsilon,\ldots,\epsilon^{m-1})$ then $ST{\mathcal H}_m(ST)^{-1}={\rm diag}(\lambda_1,\ldots,\lambda_m)+K$ with $\|K\|\leq C\epsilon$. From \eqref{eq:speccalH} we have $|{\mathsf{Im}}\,\lambda_j|\leq c_1|\epsilon|$. This proves that ${\mathsf{Re}} (iST{\mathcal H}_N(ST)^{-1}X,X)|\leq C|\epsilon||X|^2$ for any $X\in {\mathbb C}^d$. Therefore $e^{-cs\epsilon}\leq \|(ST)e^{is{\mathcal H}_N}(ST)^{-1}\|\leq e^{cs\epsilon}$ for $0\leq \epsilon\leq \epsilon_0$ with some $c>0$, $\epsilon_0>0$ because $\|STL_N(ST)^{-1}\|\leq C\epsilon$. Since $\|S^{-1}\|\leq C\epsilon^{-(m-1)}$, $\|S\|\leq C$ this proves 
\[
\epsilon^{m-1}e^{-cs\epsilon}/C
\ \leq\ \|e^{is{\mathcal H}_N}\|
\ \leq\  C\epsilon^{-(m-1)}e^{cs\epsilon}.
\]
\end{example}
\begin{example}
If $A(t,x,\xi)$ is uniformly diagonalizable then \eqref{eq:Matexp} holds with $\theta=0$. 
\rm
 Indeed, by assumption,  there exists $T=T(t,x,\xi)$ with uniform bounds of $\|T\|$ and $\|T^{-1}\|$ independent of $(t,x,\xi)$ such that $T^{-1}A(t,x,\xi)T$ is a diagonal matrix. Considering $T^{-1}e^{is{\mathcal H}_m}T$ we may assume that $i{\mathcal H}_m={\rm diag}(i\lambda_1,\ldots, i\lambda_m)+A_1$, $\|A_1\|\leq C\epsilon$ where $\lambda_j$ are real. This proves clearly 
\[
e^{-cs\epsilon}/C\leq \|e^{is{\mathcal H}_m}\|\leq Ce^{cs\epsilon}.  
\]
\end{example}

{\begin{example}
 If  for any $(t,x,\xi;\epsilon)$ there is $T=T(t,x,\xi;\epsilon)$ with uniform bounds of $\|T\|$ and $\|T^{-1}\|$ independent of $(t,x,\xi;\epsilon)$ such that $T^{-1}{\mathcal H}_mT$ is a direct sum $\Sigma\oplus A_j$ where the size of $A_j$ is at most $\mu$. Then  \eqref
{eq:Matexp} holds with $\theta=\mu-1$, \rm  which follows by a repetition of similar arguments in Example \ref{example:1}. Our results are the first results that take account of this 
purely system behavior.  That is roots of  high multiplicity but small blocks behave according to the 
size of the blocks and not the multiplicity.
\end{example}

\begin{example}\rm If there is $r\in\NN$ such that for any $(t,x,\xi,\epsilon)$ we can find $c(t,x,\xi,\epsilon)\in \CC$ such that
\[
{\rm Rank}\big({\mathcal H}_m(t,x,\xi,\epsilon)-c(t,x,\xi,\epsilon)I
\big)\leq r\,.
\]
Then hypothesis \eqref
{eq:Matexp} holds with $\theta=r$.
\end{example}

\section{Hyperbolicity and spectral bounds}

Suppose $\Omega\subset \RR^d$ is open $A(t,x)\in  C^0(\RR\,;C^{m+1}(\Omega))$ 
is an
$m\times m$ matrix 
valued function.  Assume that 
\begin{equation}
\label{eq:HypMat}
For\  all\  (t,x)\in \RR\times \Omega,
 \qquad {\rm Spectrum}\,A(t,x)\ \subset\ \RR
 \,.
\end{equation}
Define
$$
H(t,x,y,s)\
:=\
\sum_{|\alpha|\leq m}
\frac{s^{|\alpha|}}{\alpha!}\ y^{\alpha}\,
\partial_x^{\alpha}
A(t, x)\,.
$$
The values $H(t,x,y,is)$ for $y\ne 0$ and $s$ real give an extension
of $A$ to complex arguments $t,x+isy$.

The next proposition is the main result of the section giving spectral bounds
on the Taylor polynomial $H$.

\begin{proposition}
\label{pro:hypspect} Assume \eqref{eq:HypMat}. For any $T>0$ and compact set $K\subset \Omega$ there exist $\delta>0$ and $C>0$ so that
for all $x\in K$, $|t|\leq T$, $|y|\leq 1$ and $|s|\leq \delta$,
\[
\zeta\mbox{~is an eigenvalue of~}H(t,x,y,is)
\quad
\Longrightarrow
\quad
 |{\mathsf{Im}}\,\zeta|\leq C|s|
 \,.
\]
\end{proposition}

\subsection{Quantitative Nuij}

The first step in the proof of Proposition \ref{pro:hypspect}
is to prove 
 a quantitative version  of Nuij's root splitter (\cite{Nuij})
 due to Wakabayashi \cite{Waka} (see also \cite[Lemma 3.1]{CoRau}).
\begin{lemma}[Nuij]
\label{lem:nuij:anex} If $P(\zeta)$ is a monic polynomial in $\zeta$ of degree $m$ 
all of whose roots are real  define for $s\in \RR$,  real
 $\lambda_1(s)< \lambda_2(s)<\cdots < \lambda_m(s)$
so that
$(1+sd/d\zeta)^{m}P(\zeta)=\prod_{j=1}^m(\zeta-\lambda_j(s))$. 
Then there exists $c=c(m)>0$ so that for  $s\in \RR$,
\[
\lambda_{k+1}(s)-\lambda_k(s)\ \geq\
 c\,|s|,
\qquad
 j=1,\ldots,m-1.
\]
\end{lemma}
{\bf Proof:} Let $P(\zeta)=\prod_{j=1}^m(\zeta-\lambda_j)$ with $\lambda_1\leq\cdots \leq \lambda_m$ and consider for
$\ep=1,\ldots,m+1$, the successive 
Nuij splittings for $s>0$ (the case $s<0$ is similar),
\[
(1+sd/d\zeta)^{\ep-1}P(\zeta)=\prod_{j=1}^m(\zeta-\lambda_j^{\ep}(s)),
\]
where
$\lambda_1^{\ep}(s)< \cdots< \lambda_{l-1}^{\ep}(s)\leq
\lambda_l^{\ep}(s)\le\cdots\le \lambda_m^{\ep}(s)$.
Compute
\begin{equation}
\label{eq:nuij:add}
h_{\ep}(\zeta,s)
\ =\
\frac{(1+sd/d\zeta)^{\ep}P(\zeta)}{(1+sd/d\zeta)^{{\ep}-1}P(\zeta)}
\ =\
1
\ +\
s\,  \sum_{j=1}^m\,
\frac{1}{\zeta-\lambda_j^{\ep}(s)}
\ .
\end{equation}
Consider the passage from the 
roots of the denominator called mother roots to the roots of the 
numerator called daughters.
 The derivative $dh_\ep/d\zeta$ 
 is strictly negative on each interval not including
a mother
root, and 
 $\lim_{|\zeta|\to \infty} h=1$. 
 The graph of $h$ below has  four mother roots where the dotted verticals cross the 
 horizontal axis.  The mother roots toward the right may  have 
high multiplicity.
\begin{center}
\includegraphics[height=4cm, keepaspectratio]{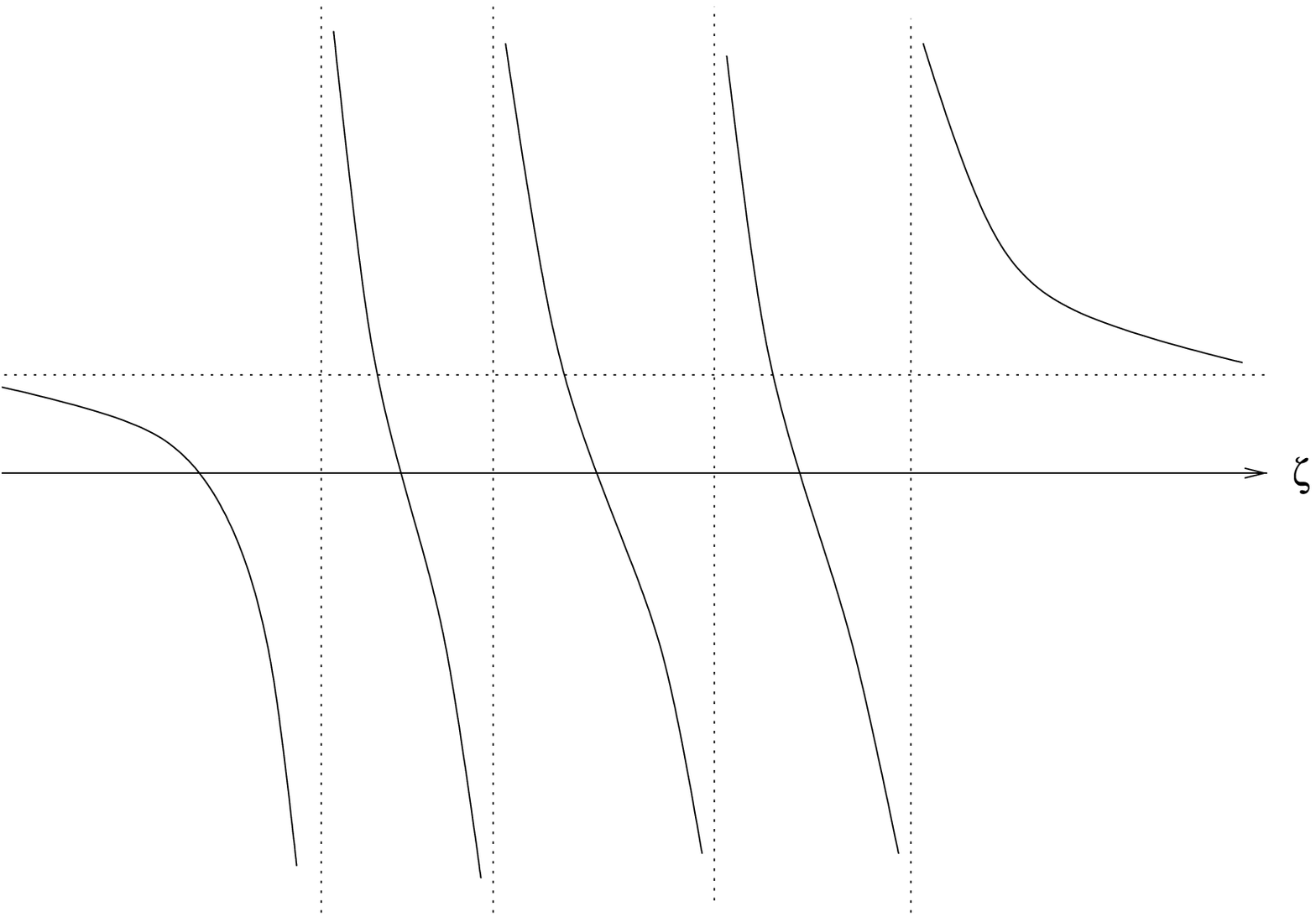}
\end{center}
There is a simple
daughter root to the 
left of the mother roots and a new simple daughter root between each of the mother roots.  
Each multiple mother root becomes a daughter root with  multiplicity reduced by  one
and gives rise to a daughter root to the left.
Each simple mother root yields a daughter to  left.  The
$\{\lambda_k^{\ep+1}(s)\}$ are all real, separate the $\{\lambda^{\ep}_k(s)\}$,
and the first ones are simple.  That is,
$$
\lambda_1^{{\ep}+1}(s)\leq \lambda_1^{\ep}(s)\leq\lambda^{{\ep}+1}_2(s)
\le \lambda^\ell_2(s)
\leq \cdots \leq \lambda^{{\ep}+1}_{m}(s)\leq \lambda_m^{\ep}(s),
$$
$$
\lambda_1^{\ep}(s)<\lambda_2^{\ep}(s)<\cdots <\lambda_{{\ep}-1}^l(s)< \lambda_{\ep}^l(s)\leq \cdots \leq \lambda_m^l(s)\,.
$$

We prove by induction on $\ep\ge 2$,
 that there exists $c_{{\ep}}>0$ such that
\begin{equation}
\label{eq:nuij:anex}
\lambda_k^{\ep}(s)\ -\
\lambda_{k-1}^{\ep}(s)\ \geq\
 c_{{\ep}}\,s,\qquad
k=2,\ldots,{\ep}.
\end{equation}
The summands $s/(\zeta - \lambda^l_j(s))$ 
in \eqref{eq:nuij:add} are all negative to the left of the mother roots.
For $l=1$ the first is equal to $-1$ when $\zeta=\lambda^1_1(s)-s$. 
Therefore
$h_1(\lambda_1^1-s,s)<0$.
 The root  $\lambda^2_1(s)$  located where the graph of $h_1$ crosses the axis
 and therefore to the left of $\lambda^1_1-s$, so
 $\lambda^2_1(s)\leq \lambda^1_1-s$.  
 
 From $\lambda^2_2(s)\geq \lambda^1_1$ it follows that
  $\lambda_2^2(s)- \lambda_1^2(s)\geq s$.
  Therefore  \eqref{eq:nuij:anex} holds with $c_2=1$ when ${\ep}=2$. 
 
 Suppose \eqref{eq:nuij:anex} holds for $2\leq k\leq {\ep}$. Prove the case
 $\ep +1$.  In  \eqref{eq:nuij:add} with $\zeta=\lambda_k^l(s)-\delta s$ the last  $m-k+1$ terms
 are negative and the first $k-1$ terms do not exceed $1/(\lambda_k^l(s)-\delta s-\lambda_{k-1}^l(s))$. Therefore by \eqref{eq:nuij:anex},
$$
h_{\ep}(\lambda_k^{\ep}(s)-\delta s,s)
\ \leq\  1+\frac{s(k-1)}{\lambda_k^{\ep}(s)-\delta s-\lambda^{\ep}_{k-1}(s)}
-\frac{1}{\delta}
\ \leq\
 1+\frac{k-1}{c_{{\ep}}-\delta}-\frac{1}{\delta}.
$$
The right hand side vanishes when $\delta=\big(k+c_{{\ep}}-\sqrt{(k+c_{{\ep}})^2-4c_{{\ep}}}\big)/2>0$. We have $h_{\ep}(\lambda_k^{\ep}(s)-\delta s,s)\leq 0$.  Therefore $\lambda_k^{{\ep}+1}(s)\leq \lambda_k^{\ep}(s)-\delta s$. Define
\[
c_{{\ep}+1}\ :=\
\min_{2\leq k\leq {\ep}} {\big(k+c_{{\ep}}-\sqrt{(k+c_{{\ep}})^2-4c_{{\ep}}}\big)/2>0}
\]
Then
\begin{align*}
\lambda_{k+1}^{{\ep}+1}(s)-\lambda_k^{{\ep}+1}(s)
\ &=\ \lambda_{k+1}^{{\ep}+1}(s)\ -\ 
\lambda_k^{\ep}(s)\ +\ 
\lambda_k^{\ep}(s)\ -\ \lambda^{{\ep}+1}_k(s)
\cr
\ &\geq\
 \lambda_k^{\ep}(s)-\lambda^{{\ep}+1}_k(s)\ \geq\
  c_{{\ep}+1}\,s
\end{align*}
for $k=1,\ldots,{\ep}$. This completes the inductive step,
so  yields \eqref{eq:nuij:anex} for ${\ep}=m+1$.
\hfill
\qed

\subsection{Three lemmas}

This subsection presents three lemmas needed in the proof of
Proposition \ref{pro:hypspect}.
  Define  
$$
Q(\zeta,t,x,y,s)  \ :=\ 
{\rm det}\,\big(\zeta I-H(t,x,y,s)\big)\,.
$$
Then  $Q(\zeta,t,x,0,s)={\rm det}\,(\zeta-A(t,x))$ and {\it for real $t,x,y,s$},
\begin{align*}
q(\zeta,t,x,y,s)
\ &=\
{\rm det}(\zeta I-A(t,x+sy))\ =\
{\rm det}(\zeta I-H+R_{m+1})\\
\ &=\ Q(\zeta,t,x,y,s)+R(\zeta,t,x,y,s)
\end{align*}
where $R(\zeta,t,x,y,s)$ is a polynomial in $\zeta$ of degree $m-1$ with coefficients $O(|s|^{m+1})$. 

The next lemma examines what happens when
the Taylor expansion and root splitter are applied simultaneously.  
Apply Nuij's root splitter to obtain polynomials with distinct roots denoted with a tilde,
\begin{equation}
\label{eq:Nkinji}
\begin{split}
{\tilde q}(\zeta,t,x,y,s):=(1+s\partial/\partial\zeta)^{m}q(\zeta,t,x,y,s)=\prod_{j=1}^m(\zeta-{\tilde\lambda}_j(t,x,y,s)),\\
{\tilde Q}(\zeta,t,x,y,s):=(1+s\partial/\partial\zeta)^{m}Q(\zeta,t,x,y,s)=\prod_{j=1}^m(\zeta-{\tilde\Lambda}_j(t,x,y,s)).
\end{split}
\end{equation}
\begin{lemma}
\label{lem:jitu} If $I\times K\subset \RR\times\Omega$ is compact,
 there is $s_0>0$ so that 
for  $(t,x,s)\in I\times K\times[-s_0,s_0]$ and $|y|\leq 1$,
all roots $\zeta$ of ${\tilde Q}=0$ are real.
\end{lemma}
{\bf Proof:} We may assume that $x+sy\in \Omega$ when $(t,x)\in I\times K$, $|y|\leq 1$ and $|s|\leq s_0$. 
The definitions \eqref{eq:Nkinji} imply that
 ${\tilde q}(\zeta,t,x,y,s)-{\tilde Q}(\zeta,t,x,y,s)={\tilde R}$ where ${\tilde R}(\zeta,t,x,y,s)$ is a polynomial in $\zeta$ of degree $m-1$ with coefficients $O(|s|^{m+1})$ uniformly in $(t,x)\in I\times K$, $|y|\leq 1$. 
Lemma \ref{lem:nuij:anex} implies
\[
|{\tilde\lambda}_{j+1}(t,x,y,s)-{\tilde\lambda}_j(t,x,y,s)|\ \geq\
 c(m)|s|\,.
\]
Let $C_j$ be the circle of radius $c(m)|s|/2$ with center ${\tilde\lambda}_j(t,x,y,s)$ so that $
|{\tilde q}(\zeta,t,x,y,s)|\geq (c(m)/2)^m|s|^m$ if $\zeta\in C_j$. 
Since $
|{\tilde q}(\zeta,t,x,y,s)-{\tilde Q}(\zeta,t,x,y,s)|\leq C|s|^{m+1}$,
 Rouch\'e's theorem implies that there exists $s_1>0$ such that  there is
 exactly one root  of ${\tilde Q}(\zeta,t,x,y,s)$ inside $C_j$ for $|s|\leq s_1$.
Since ${\tilde Q}(\zeta,t, x,y,s)$ is a real polynomial,
the root must be real.  
\hfill
\qed
\begin{lemma}
\label{lem:tQroot} 
Suppose that
 ${\tilde Q}({\bar\lambda},{\bar t},{\bar x},0,0)={\rm det}({\bar \lambda} I-A({\bar t}, {\bar x}))=0$. 
 Then there exists $\delta>0$ such that when
 $|\zeta-{\bar\lambda}|<\delta$, $|t-{\bar t}|<\delta$, $|x-{\bar x}|<\delta$, $|y|<\delta$, $|s|<\delta$,
 one has
 ${\tilde Q}(\zeta,t,x,y,s)\neq 0$ if ${\mathsf{Im}}\,\zeta\leq 0$, ${\mathsf{Im}}\,s<0$ {\rm (}or ${\mathsf{Im}}\,\zeta\geq 0$, ${\mathsf{Im}}\,s>0$).
\end{lemma}
{\bf Proof.} 
Define $p(\zeta,t,x):={\rm det}(\zeta-A(t,x))=\prod_{j=1}^m(\zeta-\lambda_j(t,x))$. If ${\mathsf{Im}}\,\zeta<0$, ${\mathsf{Im}}\,s\leq 0$ then
\[
{\tilde Q}(\zeta,t,x,0,s)\ =\ (1+s\partial/\partial\zeta)^{m}p(\zeta,t,x)\ \neq\  0.
\]
Indeed,
\[
\frac{
(1+s\partial/\partial\zeta)p(\zeta,t,x)}
{p(\zeta,t,x)}
\ =\ 1+s\sum_{k=1}^m1/(\zeta-\lambda_j(t,x))
\ =\
0
\]
implies that 
$\sum_{k=1}^m1/(\zeta-\lambda_j(t,x))=-1/s$ so that ${\mathsf{Im}}\,\sum_{k=1}^m1/(\zeta-\lambda_j(t,x))>0$ provided  that ${\mathsf{Im}}\,\lambda_j(t,x)\geq 0$ for all $j$, which is a contradiction. 

That is
 $(1+s\partial/\partial\zeta)p(\zeta,t,x)=0$ implies ${\mathsf{Im}}\,\zeta\geq 0$. It is enough to repeat this argument. Since ${\tilde Q}({\bar \lambda},{\bar t}, {\bar x}, 0,0)=0$ and ${\tilde Q}(\zeta,t,x,0,s)$ is a polynomial in $s$ of degree $m$ with leading term $ms^m$, we can find $\delta_1>0$ so that the roots $s$ of 
\[
{\tilde Q}(\zeta,t,x,y,s)\ =\ 0
\]
 with $|s|<s_0$ are continuous in $(\zeta,t,x,y)$ for $|\zeta-{\bar \lambda}|<\delta_1$, $|t-{\bar t}|<\delta_1$, $|x-{\bar x}|<\delta_1$, $|y|<\delta_1$. 
 
 Suppose that ${\tilde Q}({\hat\zeta},{\hat t}, {\hat x}, {\hat y}, {\hat s})=0$ with ${\mathsf{Im}}\,{\hat\zeta}\leq 0$, ${\mathsf{Im}}\,{\hat s}<0$, $|{\hat s}|\leq s_0$, $|{\hat \zeta}-{\bar \lambda}|<\delta_1$, $|{\hat t}-{\bar t}|<\delta_1$, $|{\hat x}-{\bar x}|<\delta_1$, $|{\hat y}|<\delta_1$.  Moving ${\hat \zeta}$ little bit
  if necessary, we may assume that ${\mathsf{Im}}\,{\hat\zeta}<0$.  Consider $F(\theta)=\min_{|s(\theta)|\leq s_0} {{\mathsf{Im}}\,s(\theta)}$ where the minimum is taken over all roots $s(\theta)$ of ${\tilde Q}({\hat\zeta},{\hat t}, {\hat x}, \theta{\hat y}, s)=0$ with $|s(\theta)|\leq s_0$. Since $F(1)<0$, $F(0)\geq 0$ there exist ${\hat \theta}$ and $s({\hat\theta})$ such that ${\mathsf{Im}}\,s({\hat\theta})=0$ which contradicts  Lemma \ref{lem:jitu}. 
  
  The proof for the case ${\mathsf{Im}}\,\zeta\geq 0$, ${\mathsf{Im}}\,s>0$ is similar.
\qed
\begin{lemma}
\label{lem:aa} Assume \eqref{eq:HypMat}. Let $({\bar t}, {\bar x})\in \RR\times \Omega$ and let ${\bar\lambda}$ be an eigenvalue of $A({\bar t}, {\bar x})$ with multiplicity $r$ so that ${\rm det}({\bar\lambda}-A({\bar t}, {\bar x}))=0$. Then there exist $\delta>0$ and $c>0$ so that 
for all $|\lambda-{\bar\lambda}|\leq \delta$, $|t-{\bar t}|<\delta$, $|x-{\bar x|\leq \delta}$, $|y|\leq \delta$ and $|s|\leq \delta$,
\begin{equation}
\label{eq:natumi}
|Q(\lambda+is,t,x,y,is)|\geq c\,|s|^r\,.
\end{equation}
\end{lemma}
{\bf Proof.} For $|t-{\bar t}|\leq \delta$, $|x-{\bar x|\leq \delta}$, $|y|\leq \delta$, $|s|<\delta$
define
$I:=\{i\mid {\tilde \Lambda}_i({\bar t}, {\bar x}, 0,0)={\bar \lambda}\}$ and $I^c:=\{i\mid {\tilde \Lambda}_i({\bar t}, {\bar x}, 0,0)\neq {\bar \lambda}\}$.  Then,
\begin{align*}
{\tilde Q}(\zeta,t,x,y,is)
\ &=\
\prod_{j\in I}(\zeta-{\tilde\Lambda}_j(t,x,y,is))
\ \prod_{j\in I^c}(\zeta-{\tilde\Lambda}_j(t,x,y,is))\\
\ &:=\
{\tilde Q}_1(\zeta,t,x,y,is)\
{\tilde Q}_2(\zeta,t,x,y,is)\,.
\end{align*}

Lemma \ref{lem:tQroot}  implies that $\pm\,{\mathsf{Im}}\,{\tilde\Lambda}_j(t,x,y,is)\geq 0$ if $\pm\,s<0$ and $j\in I$. This shows that if $M>0$, then 
\[
|{\tilde Q}_1(\lambda+iMs,t,x,y,is)|
\ \geq\
 2^{-r/2}\
 \prod_{j\in I}\big(|\lambda-{\mathsf{Re}}{\tilde\Lambda}_j|+M|s|+|{\mathsf{Im}}{\tilde\Lambda}_j|\big)
\]
for small $s\in\RR$. 
The right-hand side is  bounded from below by 
\begin{equation}
\label{eq:Qone:1}
c(M|s|)^k\
\sum 
\prod_{j_p\in I, j_1<\cdots<j_{r-k}}\big(|\lambda-{\mathsf{Re}}{\tilde\Lambda}_j|+M|s|+|{\mathsf{Im}}{\tilde\Lambda}_j|\big)
\end{equation}
 for all $1\leq k\leq r$. We prove that there are $c_k$ such that
\begin{equation}
\label{eq:QandQtilde}
Q(\zeta,t,x,y,s)={\tilde Q}(\zeta,t,x,y,s)+\sum_{l=1}^mc_l(s\partial/\partial\zeta)^l{\tilde Q}(\zeta,t,x,y,s).
\end{equation}
The definition of $\tilde Q$ implies
 $$
 (1-s\partial/\partial\zeta)^m{\tilde Q}
 \ =\
  (1-s\partial/\partial\zeta)^m (1+s\partial/\partial\zeta)^mQ
\  =\
 (1-s^2\partial^2/\partial\zeta^2)^mQ\,.
 $$
  Repeating this argument  yields 
\begin{align*}
(1+s^{2\ep}\partial^{2\ep}/\partial\zeta^{2\ep})^m\cdots(1+s^2\partial^2/\partial\zeta^2)(1-s\partial/\partial\zeta)^m{\tilde Q}
=(1-s^{4\ep}\partial^{4\ep}/\partial\zeta^{4\ep})^mQ
\end{align*}
where the right-hand side coincides with $Q$ if  $4\ep\geq m+1$. 

For $|s|M\leq 1$, note that
\begin{align*}
|\big((s\partial/\partial\zeta)^k{\tilde Q}_1\big)(\lambda&+iMs,t,x,y,is)|
\\
\ &\lesssim\
 \sum |s|^k
\prod_{ j_p\in I,\,j_1<\cdots<j_{r-k}}\big(|\lambda-{\mathsf{Re}}\,{\tilde\Lambda}_j|+M|s|+{|\mathsf{Im}}\,{\tilde\Lambda}_j|)\\
\ &\lesssim\
 M^{-k}|{\tilde Q}_1(\lambda+iMs,t,x,y,is)|
\end{align*}
by \eqref{eq:Qone:1} and
\[
|\big((s\partial/\partial\zeta)^k{\tilde Q}_2\big)(\lambda+iMs,t,x,y,is)|
\ \leq \
C|s|^k.
\]
Leibniz' rule yields
\begin{align*}
\big|(s\partial/\partial\zeta)^l\big({\tilde Q}_1{\tilde Q}_2)(\lambda&+iMs,t,x,y,is)\big|
\cr
\ &\lesssim\ 
 M^{-l}\sum_{j=0}^l (M|s|)^{l-j}
\  |{\tilde Q}_1(\lambda+iMs,t,x,y,is)|
\cr
\ &\lesssim \ M^{-l}|{\tilde Q}_1(\lambda+iMs,t,x,y,is)|
\end{align*}
because $M|s|\leq 1$. Therefore using \eqref{eq:QandQtilde}, $|Q(\lambda+iMs,t,x,y,is)|$ is bounded from below by
\begin{align*}
 |{\tilde Q}_2(\lambda&+iMs,t,x,y,is)|
\Big\{|{\tilde Q}_1(\lambda+iMs,t,x,y,is)|\\
&-C\sum_{l=1}^mM^{-l}|{\tilde Q}_1(\lambda+iMs,t,x,y,is)||{\tilde Q}_2(\lambda+iMs,t,x,y,is)|^{-1}\Big\}\ .
\end{align*}
Choosing $M>0$ large yields
\[
|Q(\lambda+iMs,t,x,y,is)|\ \geq\
 c|{\tilde Q}_1(\lambda+iMs,t,x,y,is)|
 \ \geq\ cM^r|s|^r
\]
because 
\[
|{\tilde Q}_2(\lambda+iMs,t,x,y,is)|\ =\
\prod_{j\in I^c}|\lambda+iMs-{\tilde \Lambda}_j(t,x,y,is)|
\ \geq\
 c_1
\ >\ 0.
\]
Since 
\begin{align*}
Q(\lambda+is,t,x,y,is)\ =\ 
Q(\lambda+iM(M^{-1}s),t,x,My,iM^{-1}s)
\ \geq\  c\, |s|^r
\end{align*}
the desired conclusion follows.
\hfill
\qed
\vskip.2cm
%

%
%
%
%
%
%

\subsection{Proof of Propsition \ref{pro:hypspect}  }

{\bf Proof.}   Suppose that  $({\bar t}, {\bar x})\in \{|t|\leq T\}\times K$ and ${\bar\lambda}_j$
are the distinct eigenvalues of $A({\bar t}, {\bar x})=H({\bar t}, {\bar x},0,0)$, possibly with multiplicity
greater than one. Then there is $\delta>0$ such that Lemma \ref{lem:aa} holds for any $j$. Taking $0<\delta_1\leq \delta$ small one can assume that $|{\mathsf{Re}}\,\zeta-{\bar\lambda}_{\mu}|<\delta$ 
for  some $\mu$ if $Q(\zeta,t,x,y,is)=0$ and $|t-{\bar t}|\leq \delta_1$, $|x-{\bar x}|\leq \delta_1$, $|y|\leq \delta_1$, $|s|<\delta_1$. 

Suppose that there were $|{\hat t}-{\bar t}|\leq \delta_1$, $|{\hat x}-{\bar x}|\leq \delta_1$, $|{\hat y}|\leq \delta_1$, $|{\hat s}|<\delta_1$ and $\zeta_j$ such that 
\[
\big| {\mathsf{Im}}\,\zeta_j({\hat t}, {\hat x}, {\hat y},{\hat s}) \big|
\ >\ |{\hat s}|.
\]
Clearly ${\hat y}\neq 0$ and
${\hat s}\neq 0$. First suppose that ${\mathsf{Im}}\,\zeta_j({\hat t}, {\hat x}, {\hat y},{\hat s})>|{\hat s}|$. Introduce
\begin{equation}
\label{eq:okii}
\Lambda(\theta)\ :=\
\max\Big\{{\mathsf{Im}}\,\zeta_j({\hat t}, {\hat x}, \theta{\hat y}, {\hat s})\ :\
|{\mathsf{Re}}\,\zeta_j\,-\,      {\bar \lambda}_{\mu}|<\delta
\Big\}\,.
\end{equation}
 Note that $\Lambda(0)=0$ and $\Lambda(1)>|{\hat s}|$. Since $\Lambda(\theta)$ is continuous there exist $\ep$ and  ${\hat \theta}$ such that $\Lambda({\hat \theta})=|{\hat s}|$ so that $\zeta_{\ep}({\hat t}, {\hat x}, {\hat \theta}{\hat y}, {\hat s})=\al+i|{\hat s}|$ with $\al\in\RR$ and $Q(\al+i|{\hat s}|, {\hat t}, {\hat x}, {\hat \theta}{\hat y}, i{\hat s})=0$.
 This contradicts  Lemma \ref{lem:aa} if ${\hat s}>0$.  
 
 If ${\hat s}<0$ then  $H({\hat t}, {\hat x}, {\hat \theta}{\hat y}, i{\hat s})=H({\hat t}, {\hat x}, -{\hat \theta}{\hat y}, -i{\hat s})$
 yields $Q(\al-i{\hat s}, {\hat t}, {\hat x}, -{\hat \theta}{\hat y}, -i{\hat s})=0$.
 This  contradicts  Lemma \ref{lem:aa}.
 
  If ${\mathsf{Im}}\,\zeta_j({\hat t}, {\hat x}, {\hat y}, {\hat s})<-|{\hat s}|$ it is enough to consider the minimum in \eqref{eq:okii}.
%
%
Thus we conclude that if $Q(\zeta, t,x,y,is)=0$ with $|t-{\bar t}|\leq \delta_1$, $|x-{\bar x}|\leq \delta_1$, $|y|\leq \delta_1$, $|s|<\delta_1$ then $|{\mathsf{Im}}\,\zeta|\leq |s|$. Since $\{|t|\leq T\}\times K$ is compact there is $\delta_2>0$ such that $|{\mathsf{Im}}\,\zeta|\leq |s|$ if $Q(\zeta, t,x,y,is)=0$ and $|t-{\bar t}|\leq \delta_2$, $|x-{\bar x}|\leq \delta_2$, $|y|\leq \delta_2$, $|s|<\delta_2$. 
The identity
\[
H(t,x,y,is)=H(t,x,\delta_2 y,i\delta_2^{-1}s)
\]
yields the desired
conclusion.
\hfill
\qed

\section{The symmetrizer construction}

\subsection{Lyapunov function for linear ODE}
\label{sec:lyapunov}

Suppose that $M$ is a matrix all of whose eigenvalues
lie in the open left half plane $\{{\rm Re}\,z<0\}$.
The solutions $X(t)$ of the ordinary differential equation
$$
X^\prime \ =\ M\,X
$$
tend exponentially to zero as $t\to \infty$.

Lyapunov proved that there are positive definite symmetric
matrices $R$  so that the scalar product
$(RX,X)$ is strictly decreasing on orbits.
For differential equations the quantity $(R\cdot\,,\,\cdot)$ is called
a Lyapunov function.   In the partial differential equations context,
$R$ is often called a symmetrizer.

There  is  a very 
 clever explicit choice 
\begin{equation}
\label{eq:R}
R\ =\
\int_0^\infty
(e^{sM})^*
\ e^{sM}\ ds\,.
\end{equation}
For that $R$, compute
\begin{align*}
RM +M^*R
&=\int_0^{\infty}e^{sM^*}e^{sM}M\ ds\ 
 +\int_0^{\infty}M^*e^{sM^*}e^{sM}\ ds\\
&=\int_0^{\infty}e^{sM^*}\ \frac{d}{ds}(e^{sM})\ ds\
+\int_0^{\infty}\frac{d}{ds}(e^{sM^*})\ e^{sM}\ ds\\
&=\int_0^{\infty}\frac{d}{ds}(e^{sM^*}e^{sM})\ ds
\ =\
-I.
\end{align*}
Therefore
\begin{align*}
\frac{d}{dt}
(RX(t),X(t))
 &=
(RX^\prime(t),X(t))  + 
(RX(t),X^\prime(t))\\
&= 
(RMX(t),X(t))  +
(RX(t),MX(t))
\\
&= 
\Big(
\big(
RM +M^*R\big)X(t)
\,,\,
X(t)
\Big)
\ =\ -\,\big(X(t),X(t)\big),
\end{align*}
proving that $(R\, \cdot\, ,  \,  \cdot)$ is a 
strict
Lyapunov function.

The last identity is easily understood.   With 
$X(t)=e^{tM}X(0)$,   the definition of $R$ yields the identities
for $s>0$.
$$
\big(RX(0),X(0) \big) 
=
\int_0^\infty \|X(t)\|^2\ dt\,,
\qquad
\big(R X(s), X(s)\big) 
=
\int_s^\infty \|X(t)\|^2\ dt\,,
$$
and the formula for $(RX,X)^\prime$ follows.

For applications to partial differential equations one has
matrices $M$ that depend smoothly on parameters and it is
important that the symmetrizers also depend smoothly.
The standard constructions of lyapunov functions
depending either on Schur's unitary upper triangularization
or Jordan's canonical upper triangularization do not have
smooth dependence.  
Formula \eqref{eq:R}
in contrast does depend smoothly on parameters. It
pays no attention  
to the spectral details of $M$.  
Where eigenvalues  cross
and the associated spectral projections usually misbehave,
the  formula for $R$ does not.

The  identity $RM+M^*R<0$  is important.
It implies a negativity of symbols that translates, thanks to 
the sharp
G\aa rding inequality,
 to a negativity of operators in our 
 application.

\subsection{Symmetrizer $R$ and its derivatives }

Assume \eqref{eq:Matexp} and hence \eqref{eq:expH}. 
Define
\[
M(a,\ell,\tau,\rho, t,x,\xi)\ :=\
iH_N(\ell,\tau,t,x,\xi)
\ -\
a\, \lr{\xi}_\ell^{\rho}
\]
with  
\begin{equation}
\label{eq:contraint5}
0<\rho<1\le \min\{a,\ell\}
\end{equation}
Proposition \ref{pro:hypspect} implies that there is an $a_0\ge 1$, $c>0$ so that 
$$
{\rm Spectrum} \, M
\ \subset\
\Big\{z\ :\ {\rm Re}\,z \,\le\, c(a_0-a)\,\langle\xi\rangle^\rho_\ell
\Big\}\,.
$$
We suppose that 
\begin{equation}
\label{eq:constraint8}
a\ \ge\ a_0+1\,.
\end{equation}
The parameters  $\tau, a,T$ are constrained to satisfy
\begin{equation}
\label{eq:rangetau}
c_1
\ \leq\
 c\,\tau
 \ \leq\
  T,\qquad  2\ c\ T\leq a
\end{equation}
for some $T>c_1>0$ .
For ease of reading, the $\ell,\tau,a,\rho$ dependence of 
$M$ and $R$ is often omitted. 
Introduce the candidate symmetrizer
$$
R(a,\ell,\tau,\rho, t,x,\xi)
\ :=\ 
a\int_0^\infty
\langle\xi\rangle^\rho_\ell\
\big(
e^{sM(t,x,\xi)}\big)^*
\big(
e^{sM(t,x,\xi)}\big)
\ ds.
$$
We need lower bounds on $R$ so that it yields good estimates
and need to verify that  $R$ defines
a classical  symbol.   Interestingly, we do not need that $R$ is a Gevrey symbol.

The parameters $\ell,\rho, a$ are constrained by 
\begin{equation}
\label{eq:aellconstraint}
1\ \leq\
 a
 \ \leq\
  \ell^{1-\rho}.
\end{equation}
%


Since $\|e^{sM}\|=e^{-as\lr{\xi}_\ell^{\rho}}\|e^{isH_N}\|$,  \eqref{eq:rangetau}
implies 
\[
\tau^{\theta}\lr{\xi}_\ell^{-\theta(1-\rho)}e^{-c_1as\lr{\xi}_\ell^{\rho}}/C
\ \leq\
 \big\|e^{sM}\big\|
 \ \leq\
  C \tau^{-\theta}\lr{\xi}_\ell^{\theta(1-\rho)}\ e^{-c_2a s\lr{\xi}_\ell^{\rho}}
\]
with  $c_i,C>0$ independent of $\ell,\tau,a,t,x,\xi,s$. 
This yields 
\begin{align*}
(Rv,v)\ & =\
a  \int_0^{\infty}\lr{\xi}_h^{\rho}\|e^{sM}v\|^2\,ds\\
\ &\geq\  C^{-2}\tau^{2\theta}\|v\|^2\lr{\xi}_{\ell}^{-2\theta(1-\rho)}
 \int_0^{\infty}a\lr{\xi}_\ell^{\rho}e^{-2c_1as\lr{\xi}_\ell^{\rho}}\,ds\\
\ &\geq\
 c'\,\tau^{2\theta}
\lr{\xi}_\ell^{-2\theta(1-\rho)}\|v\|^2\,.
\end{align*}
This is equivalent to the important lower bound
\begin{equation}
\label{eq:lowerboundR}
R
\ \geq\
 c'\,\tau^{2\theta}\lr{\xi}_{\ell}^{-2\theta(1-\rho)}.
\end{equation}
\begin{theorem}
\label{thm:symbol}
Assume  \eqref{eq:Matexp}, \eqref{eq:rangetau} with $K=\RR^d$ with
$0\leq \theta\le m-1$.  Denote 
$\nu:=\theta(1-\rho)$.   Suppose 
that $A(t,x,\xi)$ is lipschitzian in time uniformly on compact sets with values in the ${\tilde S}^1(\RR^d\times\RR^d)$. Then 
$R(t,x,\xi)$
(resp. $\partial_tR$)  is bounded in time uniformly on compacts with values in 
${\tilde S}^{2\nu}_{\rho-\nu,1-\rho+\nu}(\RR^d\times\RR^d)$ and 
(resp. ${\tilde S}^{1-\rho+3\nu}_{\rho-\nu,1-\rho+\nu}(\RR^d\times\RR^d)$).
   That is for all $\alpha,\beta$,
\begin{equation}
\label{eq:estRR}
\begin{split}
&|\partial_x^{\beta}\partial_{\xi}^{\alpha}R(t,x,\xi)|
\ \leq\
 C_{\alpha\beta}a^{-|\alpha+\beta|}\lr{\xi}_\ell^{2\nu+(1-\rho+\nu)|\beta|-(\rho-\nu)|\alpha|},
\\
&|\partial_x^{\beta}\partial_{\xi}^{\alpha}\partial_tR(t,x,\xi)|
\ \leq\
 C_{\alpha\beta}a^{-|\alpha+\beta|-1}
\lr{\xi}_\ell^{1-\rho+3\nu+(1-\rho+\nu)|\beta|-(\rho-\nu)|\alpha|}
\end{split}
\end{equation}
with  $C_{\al\beta}$ independent of $a,\rho,\ell,\tau, t,x,\xi$. 
\end{theorem}
%
%
%
%

\begin{remark}
\label{rem:DtR}
\rm
The estimate for $\partial_x^{\beta}\partial_{\xi}^{\alpha}\partial_tR$ is exactly the same as the estimate
for an a derivative $\partial_x^{\gamma}\partial_{\xi}^{\alpha}R$ with $|\gamma| =|\beta|+1$.
The time derivative is like an extra space derivative.
\end{remark}

{\bf Proof.} Denote 
$$
X(s;t,x,\xi)\ :=\
e^{sM(t,x, \xi)}v\,,
\qquad
X^{\alpha}_{\beta}(s;t,x,\xi)\ :=\
\partial_x^{\beta}\partial_{\xi}^{\alpha}X(s;t,x,\xi)
\,.
$$

\vskip.1cm
{\bf Step I.  Estimates for $X^\alpha_\beta$.}    We prove, by induction on $|\alpha+ \beta|$, that
\begin{equation}
\label{eq:estX}
\big|X^{\alpha}_{\beta}(s)  \big|
\ \leq\  C_{\alpha\beta}
\big(s+\lr{\xi}_{\ell}^{-1}\big)^{|\alpha|}\ 
\big(1+s\lr{\xi}_{\ell}\big)^{|\beta|}\
\lr{\xi}_{\ell}^{\nu(|\alpha+\beta|+1)}\ 
e^{-cas\lr{\xi}_{\ell}^{\rho}}\,.
\end{equation}
The constraint \eqref{eq:aellconstraint} implies that
\begin{equation}
\label{eq:constraint9}
a\lr{\xi}_\ell^{\rho-1}\leq 1,\quad 
{\rm and}
\quad
1\leq a^{-1}\lr{\xi}_\ell^{1-\rho}.
\end{equation}
The identity
$\lr{\xi}_{\ell}=\ell\lr{\xi/\ell}$ from \eqref{eq:defanglexi} implies 
\[
|\partial_{\xi}^{\alpha}\lr{\xi}_\ell^s|
\ \lesssim\
 \lr{\xi}_\ell^{s-|\alpha|}
\qquad
{\rm so},
\qquad
\partial_{\xi}^{\alpha}\lr{\xi}_{\ell}\ =\
\ell\
\ell^{-|\alpha|}\
(\partial/\partial\zeta)^{\alpha}\lr{\zeta}\big|_{\zeta=\xi/\ell}
\,.
\]

Introduce
\[
E(s)\ :=\
\lr{\xi}_{\ell}^{\nu}\ e^{-cs\,a\lr{\xi}_{\ell}^{\rho} }
\]
so that $|X|\lesssim E(s)$ and $E(s)E({\tilde s})=\lr{\xi}_{\ell}^{\nu}E(s+{\tilde s})$.
The desired estimate 
\eqref{eq:estX} is equivalent to
\begin{equation}
\label{eq:kigen}
|X^{\alpha}_{\beta}(s)|
\ \leq\
 C_{\alpha\beta}\,
 \big(s+\lr{\xi}_{\ell}^{-1}\big)^{|\alpha|}
 \ 
 \big(1+s\lr{\xi}_{\ell}\big)^{|\beta|}\\
\
 \lr{\xi}_{\ell}^{\nu|\alpha+\beta|}\
 E(s)\,.
\end{equation}
Equations \eqref{eq:rangetau} and \eqref{eq:constraint9} imply
\begin{equation}
\label{eq:Mhyoka}
|M^{(\al)}_{(\beta)}|
\ \leq\ 
 C_{\al\beta}\lr{\xi}_{\ell}^{1-|\al|}\,.
\end{equation}

For  $|\alpha|=1$ differentiate the equation for $X$ to find,
\begin{equation}
\label{eq:Xdot}
{\dot X}^{\alpha}\ =\ 
MX^{\alpha}+M^{(\alpha)}X\,,
\qquad
X^\alpha(0)\ =\  0\,.
\end{equation}
 Then  \eqref{eq:Mhyoka} and Duhamel's representation yield
\begin{align*}
|X^{\alpha}(s)|
\ &=\
\Big|\int_0^se^{(s-{\tilde s})M}M^{(\alpha)}X\ 
d{\tilde s}\Big|
\ \lesssim\  \int_0^sE(s-{\tilde s})\,E({\tilde s})\ 
d{\tilde s}
\\
\ &=\ 
s\, \lr{\xi}_{\ell}^{\nu}\ E(s)
\ \leq\  \big(s+\lr{\xi}_{\ell}^{-1} \big)\  \lr{\xi}_{\ell}^{\nu}  \ E(s).
\end{align*}
Similarly for  $|\beta|=1$,
 ${\dot X}_{\beta}=MX_{\beta}+M_{(\beta)}X$ with $X_\beta(0)=0$ so 
\begin{align*}
|X_{\beta}|
\ \leq\
 \int_0^s E(s-{\tilde s})\lr{\xi}_{{\ell}}E({\tilde s})d{\tilde s}\ \leq\  s\lr{\xi}_{\ell}\lr{\xi}_{\ell}^{\nu}E(s)
\ \leq\  (1+s\lr{\xi}_{\ell})\lr{\xi}_{\ell}^{\nu}E(s).
\end{align*}
This proves \eqref{eq:kigen} for $|\alpha +\beta|=1$.

Assume  $k\ge 1$ and that \eqref{eq:kigen}
holds for $|\alpha+\beta|\le k$.  It suffices to prove \eqref{eq:kigen} 
$X^\gamma_\delta$ with
 $|\gamma+\delta|=k+1$. 
 Differentiation of the equation for $X$ yields
\begin{equation}
\label{eq:dotXab}
{\dot X}^{\alpha}_{\beta}
\ =\ 
M\,X^{\alpha}_{\beta}\
+\sum_{ 
{\alpha_1+\alpha_2=\alpha\,,\,
\beta_1+\beta_2=\beta}
\atop
{\alpha_1+\beta_1\ne 0}
}
 C_{\alpha_1,\beta_1}\
M^{(\alpha_1)}_{(\beta_1)}\ X^{\alpha_2}_{\beta_2}\,,
\qquad
X^\alpha_\beta(0)=0
\,.
\end{equation}
Duhamel's formula yields
$$
\big|X^{\alpha}_{\beta}(s)\big|\ \lesssim\  
\sum_{\al_1+\beta_1\ne 0}
\int_0^s \big|e^{(s-{\tilde s})M}\ 
M^{(\alpha_1)}_{(\beta_1)}\ 
X^{\alpha_2}_{\beta_2}
\big|\ d{\tilde s}\,.
$$
The inductive hypothesis estimates the right hand side by
\begin{equation}
\label{eq:extMX}
\begin{split}
&\lesssim
\sum_{\al_1+\beta_1\ne 0}\int_0^s E(s-{\tilde s})\lr{\xi}_{\ell}^{1-|\alpha_1|}({\tilde s}+\lr{\xi}_{\ell}^{-1})^{|\alpha_2|}(1+{\tilde s}\lr{\xi}_{\ell})^{|\beta_2|}
\lr{\xi}_{\ell}^{\nu|\alpha_2+\beta_2|}E({\tilde s})d{\tilde s}
\\
&\lesssim  \sum_{\al_1+\beta_1\ne 0}
s\,
\lr{\xi}_{\ell}^{1-|\alpha_1|}\
(s+\lr{\xi}_{\ell}^{-1})^{|\alpha_2|}\
(1+s\lr{\xi}_{\ell})^{|\beta_2|}\
 \lr{\xi}_{\ell}^{\nu(|\alpha_2+\beta_2|+1)}\
 E(s)\,.
\end{split}
\end{equation}
If $|\beta_1|\geq 1 $ so that $|\beta_2|<|\beta|$, then
\begin{align*}
s\ \lr{\xi}_{\ell}\ \lr{\xi}_{\ell}^{\nu(|\al_2+\beta_2|+1)}\
(1+s\lr{\xi}_{\ell})^{|\beta_2|}
\ \leq\
 \lr{\xi}_{\ell}^{\nu|\al+\beta|}\
 (1+s\lr{\xi}_{\ell})^{|\beta|}
\end{align*}
and the right-hand side of \eqref{eq:extMX} is bounded by \eqref{eq:kigen}. If $|\beta_1|=0$ so that $\beta_2=\beta$ and $|\al_1|\geq 1$,
\begin{equation}
\label{eq:beta1zero}
\begin{aligned}
s\
\lr{\xi}_{\ell}^{1-|\al_1|}&\
(s+\lr{\xi}_{\ell}^{-1})^{|\al_2|}\
\lr{\xi}_{\ell}^{\nu(|\al_2+\beta|+1)}
\\
&\leq\
\lr{\xi}_{\ell}^{-\nu(|\al_1|-1)}\
\lr{\xi}_{\ell}^{\nu|\al+\beta|}\
s\,\lr{\xi}_{\ell}^{-(|\al_1|-1)}\
(s+\lr{\xi}_{\ell}^{-1})^{|\al_2|}\\
&\leq\
 \lr{\xi}_{\ell}^{\nu|\al+\beta|}\
 (s+\lr{\xi}_{\ell}^{-1})^{|\al|}
\end{aligned}
\end{equation}
implying  the same conclusion.  This completes the inductive proof of
\eqref{eq:kigen}

\vskip.1cm
{\bf Step II.  Estimates for $\partial_t X^\alpha_\beta$.}   
Differentiating the equation $\dot X = M X$ with respect to $t$ or with respect to $x$
are entirely parallel.  With the exception that one can  only take one temporal derivative
because $M$
is only lipschitzian in $t$.   The result is a bound for $\dot X^\alpha_\beta$ that is the 
same as the bound for $X$ with one more $x$ derivative that is 
\begin{equation}
\label{eq:dotXbound}
\big|
\dot X^\alpha_\beta
\big|
\ \lesssim\
 \big(s+\lr{\xi}_{\ell}^{-1}\big)^{|\alpha|}
 \ 
 \big(1+s\lr{\xi}_{\ell}\big)^{|\beta|+1}\\
\
 \lr{\xi}_{\ell}^{\nu(|\alpha+\beta|+1)}\
 E(s)\,.
 \end{equation}

\vskip.1cm
{\bf Step III.  Estimates for $R^\alpha_\beta$.}
Begin with the estimate
from Leibniz' rule,
\begin{align}
\label{eq:Rderiv}
|\partial_x^{\beta}\partial_{\xi}^{\alpha}R|
\lesssim \sum_{\beta_1+\beta_2=\beta\atop \alpha_1+\alpha_2+\alpha_3=\alpha} \int_0^{\infty}a\Big|
\big(\lr{\xi}_{\ell}^{\rho}\big)^{(\alpha_1)}\big(e^{sM^*}\big)^{(\alpha_2)}_{(\beta_1)}\big(e^{sM}\big)^{(\alpha_3)}_{(\beta_2)}\Big|
\ ds\,.
\end{align}
Thanks to \eqref{eq:estX}, the integrand in \eqref{eq:Rderiv} is $\,\lesssim$
\begin{align*}
a\lr{\xi}_{\ell}^{\rho-|\alpha_1|}(s+\lr{\xi}_{\ell}^{-1})^{|\alpha_2+\alpha_3|}
 (1+s\lr{\xi}_{\ell})^{|\beta_1+\beta_2|} \lr{\xi}_{\ell}^{\nu(|\beta+\alpha-\alpha_1|+2)}
e^{-cas\lr{\xi}_{\ell}^{\rho}}e^{-cas\lr{\xi}_{\ell}^{\rho}}\,.
\end{align*}
Using the pair of estimates
\begin{align*}
s+\lr{\xi}_{\ell}^{-1}
\ &=\
\big(as\lr{\xi}_{\ell}^{\rho}+a\lr{\xi}_{\ell}^{-1+\rho}\big)a^{-1}\lr{\xi}_{\ell}^{-\rho}
\ \leq\ 
 \big(as\lr{\xi}_{\ell}^{\rho}+1\big)a^{-1}\lr{\xi}_{\ell}^{-\rho},
\cr
1+s\lr{\xi}_{\ell}\ &=\
as\lr{\xi}_{\ell}^{\rho}(a^{-1}\lr{\xi}_{\ell}^{1-\rho})+1
\ \leq\ 
\big(as\lr{\xi}_{\ell}^{\rho}+1\big)a^{-1}\lr{\xi}_{\ell}^{1-\rho},
\end{align*}
yields
\begin{align*}
\big|\partial_x^{\beta}\partial_{\xi}^{\alpha}R\big|
\ \lesssim \
a^{-|\beta+\alpha|}\
\lr{\xi}_{\ell}^{q}
 \int_0^{\infty}a\,(1+as\lr{\xi}_{\ell}^{\rho})^{|\beta+\alpha|}\
 \lr{\xi}_{\ell}^{\rho}\
 e^{-2cas\lr{\xi}_{\ell}^{\rho}}\ ds\,,
\end{align*}
with
$$
q:=
(1-\rho)|\beta|-\rho|\alpha|+\nu(|\beta+\alpha|+2)
=
2\nu+(1-\rho+\nu)|\beta|-(\rho-\nu)|\alpha|\,.
$$
Use $(1+as\lr{\xi}_{\ell}^{\rho})^{|\beta+\alpha|}\,e^{-cas\lr{\xi}_{\ell}^{\rho}}\lesssim 1$ to 
find 
\begin{align*}
|\partial_x^{\beta}\partial_{\xi}^{\alpha}R|
\ \lesssim\
 a^{-|\beta+\alpha|}\
 \lr{\xi}_{\ell}^{q}\
\int_0^{\infty}a\,\lr{\xi}_{\ell}^{\rho}\ e^{-cas\lr{\xi}_{\ell}^{\rho}}\ ds
\ \lesssim\
 a^{-|\beta+\alpha|}\
 \lr{\xi}_{\ell}^{q}\,.
\end{align*}
This is the first estimate of \eqref{eq:estRR}.   

{\bf Step IV.  Estimates for $\partial_t R^\alpha_\beta$.}
As in Remark   \ref{rem:DtR}, the estimate for the time derivative
is the same as taking one additional space derivative.
 The details
are left to the reader.
 
 This completes the proof of Theorem \ref{thm:symbol}.
\hfill
\qed

\section{Theorem \ref{thm:bronshtein1}, examples and proof }
\label{Sec:shomei}

Begin with some examples illustrating the conclusion.

\begin{example}\rm
If $A(t,x,\xi)=\sum_{j=1}^dA_j(t,x)\xi_j$ is uniformly diagonalizable then one can take $\theta=0$ so that Theorem \ref{thm:bronshtein1} holds with $1<s\leq 2$. 
This is the sharp index.
In \cite{Kaji} Kajitani has proved that the Cauchy problem for uniformly diagonalizable system is well posed in ${\mathcal G}^{2}(\RR^d)$ when the coefficients are  smooth 
enough in time. For $2\times 2$ uniformly diagonalizable systems with coefficients depending only on $t$, more detailed discussions on the regularity in $t$ are found in \cite{CN}.
\end{example}

\begin{example}
\rm
We can always take $\theta=m-1$ and Theorem \ref{thm:bronshtein1} holds with $1<s\leq(4m-1)/(4m-2)$. 
\end{example}

\begin{example}
\rm
If \eqref{eq:Matexp} holds with $\theta=\mu-1$, $2\leq \mu\leq m-1$ one can choose $1<s\leq(4\mu-1)/(4\mu-2)$ in Theorem \ref{thm:bronshtein1}.   This is not sharp for $m=2$.
\end{example}

{\bf Proof of Theorem  \ref{thm:bronshtein1}.} {\bf Step I.  Compact support in $x$.}
Choose $R>0$ so that the support of $f$, $g$, and $\nabla_{t,x} A_j, \nabla_{t,x} B$
are all contained in $\{|x|\le R\}$.   Denote by $c_{max}$ an upper bound
for the propagation speed for the constant coefficient hyperbolic
operator $L$ on $|x|\ge R$.

Finite speed result applied in $|x|\ge R$ implies that $u$ vanishes
for $|x|\ge R + c_{max}t$ for $t\ge 0$.  

\vskip.1cm
{\bf Step II.  First {\it a priori} estimate.}
Consider \eqref{eq:Cauchy}.
Set $v=e^{\lr{D}_{\ell}^{\rho}(T-a t)}u$ with small $T$  to be chosen below.
Define
\[
{\tilde A}
\ :=\ e^{\lr{D}_{\ell}^{\rho}(T-a t)}Ae^{-\lr{D}_{\ell}^{\rho}(T-a t)},\qquad{\tilde B}
\ :=\
e^{\lr{D}_{\ell}^{\rho}(T-a t)}Be^{-\lr{D}_{\ell}^{\rho}(T-a t)}
\]
 and ${\tilde  f}:=e^{\lr{D}_{\ell}^{\rho}(T-a t)}f$. 
Compute
\begin{equation}
\label{eq:energy}
\begin{split}
\frac{d}{dt}\big(R\    & e^{\lr{D}_{\ell}^{\rho}(T-a t)}u,e^{\lr{D}_{\ell}^{\rho}(T-a t)}u\big)=(\partial_tR\ v,v)
+(R(i{\tilde A}+{\tilde B}-a\lr{D}_{\ell}^{\rho})v,v)
\\
&+(Rv,(i{\tilde A}+{\tilde B}-a\lr{D}_{\ell}^{\rho})v)
+(R{\tilde f}, v)+(Rv,{\tilde f})\,.
\end{split}
\end{equation}
For any $0<c_1<T$, 
$c_1\leq T-at\leq T$ for $0\leq t\leq (T-c_1)/a$.
If $T>0$ is small,
Proposition \ref{pro:weyl:1} implies that ${\tilde A}=H_N+K$ with
 \begin{equation}
 \label{eq:torii}
 K\ \in\
  {\tilde S}^{\,\max{\{\rho-N(1-\rho)\,,\,-1+\rho\}}} 
  \,.
  \end{equation}
 Corollary \ref{cor:esurei} implies
${\tilde B}\in {\tilde S}^0$. 

Since $i{\tilde A}-a\lr{D}_{{\ell}}^{\rho}=M+iK$ the 
right-hand side of \eqref{eq:energy} is
equal to
\begin{align*}
(\partial_tR\, v,v)&+\big((RM+M^*R)v, v\big)
+(R{\tilde f}, v)
\\
&+\big((R(iK+{\tilde B})+(iK+{\tilde B})^*R)v, v\big)+(Rv,{\tilde f}).
\end{align*}

Recall  that $M\in {\tilde S}^1$ and
 $R$ satisfies \eqref{eq:estRR}. Therefore
\begin{align*}
R\#M
\ +\
M^*\#R
\ =\
RM\ +\
M^*R\ +\
K_1\ =\
-\,a\lr{\xi}_{\ell}^{\rho}
\ +\ K_1
\end{align*}
where $aK_1\in {\tilde S}^{1-\rho+3\nu}_{\rho-\nu\,,\,1-\rho+\nu}$ 
with bound independent of $a$ large
(see e.g. \cite[Proposition 18.5.7]{Hobook}).  Choose
$a_1\ge a_0$ so that if $a>a_1$ one has 
$
Ca^{-1}
\ \leq\ 3a/4
$.
Then, 
\begin{align*}
-a(\lr{D}_{{\ell}}^{\rho}u,u)+{\mathsf{Re}}\,(K_1u,u)
\ &\leq\  -a\|\lr{D}_{{\ell}}^{\rho/2}u\|^2
+Ca^{-1}\|\lr{D}_{{\ell}}^{(1-\rho+3\nu)/2}u\|^2
\cr
\ &\leq\
 -(a/4)\|\lr{D}_{{\ell}}^{\rho/2}u\|^2
\end{align*}
provided 
\[
\rho\geq 1-\rho+3\nu\quad
{\rm equivalently},
\quad
    \rho\geq (1+3\theta)/(2+3\theta).
\]
Note that $a\,\partial_tR\in {\tilde S}^{1-\rho+3\nu}_{\rho-\nu\,,\,1-\rho+\nu}$ with $a$-independent bound
so 
\[
{\mathsf{Re}}\,\big(\partial_tRu,u \big)
\ \leq\ Ca^{-1}\|\lr{D}_{{\ell}}^{\rho/2}u\|^2
\]
if $2\rho\geq 1+3\nu$, that is $\rho\geq (1+3\theta)/(2+3\theta)$. 

 Using
  \eqref{eq:torii},
  $R\in {\tilde S}^{2\nu}_{\rho-\nu,1-\rho+\nu}$,
 and, ${\tilde B}\in {\tilde S}^0$  yields the pair of estimates
\begin{align*}
\big| \big((R{\tilde B}+{\tilde B}^*R)v, v\big)\big| 
\ &\leq\ C\|\lr{D}_{\ell}^{\nu}v\|^2
\ \leq\
 C{\ell}^{-(\rho-2\nu)}\|\lr{D}_{\ell}^{\rho/2}v\|^2,\\
\big| \big(i(RK-K^*R)v, v\big)\big|
\ &\leq\ C\|\lr{D}_{\ell}^{\nu+\rho/2-N(1-\rho)/2}v\|^2
\ \leq\  C'\|\lr{D}_{\ell}^{\rho/2}v\|^2
\end{align*}
because $2\nu+\rho-N(1-\rho)=(2\theta-N)(1-\rho)+\rho\leq \rho$ and $2\nu-1+\rho\leq \rho$ 
when $2\rho\geq 1+3\nu$.  In addition,
\begin{align*}
|(R{\tilde f}v,v)|+|(Rv,{\tilde f})|
\ \leq\
 2\|\lr{D}_{{\ell}}^{-3\nu}Rv\|\|\lr{D}_{{\ell}}^{3\nu}{\tilde f}\|
\ \leq\
 C\|\lr{D}_{{\ell}}^{-\nu}v\|\|\lr{D}_{{\ell}}^{3\nu}{\tilde f}\|.
\end{align*}
Thus there exist $c,C>0$ so that
\begin{equation}
\label{eq:eneid}
\frac{d}{dt}(R\, v, v)\ +\ ca\|\lr{D}_{{\ell}}^{\rho/2}v\|^2
\ \leq\  C\|\lr{D}_{{\ell}}^{-\nu}v\|\|\lr{D}_{{\ell}}^{3\nu}{\tilde f}\|\,.
\end{equation}

The definition of $R$ together with 
 \eqref{eq:lowerboundR} and  \eqref{eq:rangetau}  show that 
 if $T_1<T$ and $0\leq t\leq T_1/a$, then
$R=R^*\geq c\,\lr{\xi}_{{\ell}}^{-2\nu}$.

Introduce the metric
\begin{equation}
\label{eq:metricG}
G\ :=\
a^{-2}\Big(\lr{\xi}_{\ell}^{2(1-\rho+\nu)}|dx|^2
\ +\
\lr{\xi}_{\ell}^{-2(\rho-\nu)}|d\xi|^2\Big)\,.
\end{equation}
 Then $G/G^{\sigma}=a^{-4}\lr{\xi}_{\ell}^{2(1-2\rho+2\nu)}$. Since $H=a^{2}\lr{\xi}_{{\ell}}^{2\rho-1-4\nu}(R-c\,\lr{\xi}_{{\ell}}^{-2\nu})\in S((G/G^{\sigma})^{-1/2},G)$ and $H\geq 0$, 
 the sharp G\aa rding inequality ( \cite[Theorem 18.6.7]{Hobook}) 
 yields
\[
(Hv,v)
\ \geq\  -\,C\|v\|^2.
\]
Write
\[
H=a^{2}\lr{\xi}_{{\ell}}^{\rho-1/2-2\nu}\#(R-c\,\lr{\xi}_{{\ell}}^{-2\nu})\#\lr{\xi}_{{\ell}}^{\rho-1/2-2\nu}+K
\]
where $K\in S(1,G)$.  Introduce $u:=\lr{D}_{{\ell}}^{\rho-1/2-2\nu}v$ to find
\begin{align*}
(Hv,v)\ & =\ a^{2}((R-c\,\lr{D}_{{\ell}}^{-2\nu})u,u)
+(Kv,v)
\\
&\geq
\  -C\|v\|^2
\ =\
-\,  C\|\lr{D}_{{\ell}}^{-\rho+1/2+2\nu}u\|^2\,.
\end{align*}
Since $|(Kv,v)|\leq C\|v\|^2=C\|\lr{D}_{{\ell}}^{-\rho+1/2+2\nu}u\|^2$ it follows that
\begin{equation}
\label{eq:Rsita}
(Ru,u)-c\,\|\lr{D}_{{\ell}}^{-\nu}u\|^2
\ \geq\
 -\,Ca^{-2}\|\lr{D}_{{\ell}}^{-\rho+1/2+2\nu}u\|^2.
\end{equation}
If $-\nu\geq 2\nu+1/2-\rho$, that is 
\[
\rho
\ \geq\
 \frac{1+6\theta}
{2+6\theta}
\]
then
there is a  $c'>0,\ \ell_0>0$, so that for ${\ell}\geq {\ell}_0$. 
\[
(Ru,u)
\ \geq\
 c'\,\|\lr{D}_{{\ell}}^{-\nu}u\|^2\,.
\]
Integrating \eqref{eq:eneid} 
yields
\[
\|\lr{D}_{{\ell}}^{-\nu}v(t)\|^2
\ \leq\
 c\|\lr{D}_{{\ell}}^{\nu}v(0)\|^2
 \ +\
 2CM\int_0^t\|\lr{D}_{{\ell}}^{3\nu}{\tilde f}\|\ d\tau\,,
\]
where $M:=\sup_{0\leq \tau\leq t}\|\lr{D}_{{\ell}}^{-\nu}v(\tau)\|$. Therefore
\[
\Big(M-C\int_0^t\|\lr{D}_{{\ell}}^{3\nu}{\tilde f}\|d\tau\Big)^2
\ \leq\
 \Big(\sqrt{c}\|\lr{D}_{{\ell}}^{\nu}v(0)\|+C\int_0^t\|\lr{D}_{{\ell}}^{3\nu}{\tilde f}\|d\tau\Big)^2
\]
which gives
\[
\|\lr{D}_{{\ell}}^{-\nu}v(t)\|
\ \leq\  2\sqrt{c}\|\lr{D}_{{\ell}}^{\nu}v(0)\|
\ +\ 2C\int_0^t\|\lr{D}_{{\ell}}^{3\nu}{\tilde f}\|d\tau.
\]
This proves the following important {\it a priori} estimate.

\begin{proposition}
\label{pro:itoyu:a} 
If $\rho\geq (1+6\theta)/(2+6\theta)$
then there exist $T>0$, $a>0$ and ${\ell}_0>0$ such that for any $T_1<T$ one can find  $C>0$ such that for all $u$ so that
 $e^{\lr{D}^{\rho}_{\ell}(T-at)}
\partial_{t ,x}^\gamma
u\in L^1([0,T]\,;\,H^{3\nu}(\RR^d) )$ for $|\gamma|\le 1$,
\begin{align*}
\|\lr{D}_{{\ell}}^{-\nu}e^{\lr{D}_{\ell}^{\rho}(T-a t)}u\|
\leq C\|\lr{D}_{{\ell}}^{\nu}e^{T\lr{D}_{\ell}^{\rho}}u(0)\|
+C\int_0^t\|\lr{D}_{{\ell}}^{3\nu}e^{\lr{D}_{\ell}^{\rho}(T-a \tau)}Lu\|\,d\tau
\end{align*}
for  $0\leq t\leq T_1/a$ and ${\ell}\geq {\ell}_0$, where $\nu=\theta(1-\rho)$.
\end{proposition}

\vskip.1cm
{\bf Step III.  Second {\it a priori} estimate.}
For some values of 
$\rho$ and $\theta$,
one can improve the estimate 
for  the left hand side 
$\|\lr{D}_{{\ell}}^{-\nu}e^{\lr{D}_{\ell}^{\rho}(T-a t)}u\|$ in Proposition \ref{pro:itoyu:a}.
   Recall 
  that $\partial_tu=iA(t,x,D)u+B(t,x)u+f$ and $v=e^{\lr{D}_{\ell}^{\rho}(T-a t)}u$. Then,
%
%
\begin{equation}
\label{eq:energy:b}
\begin{split}
\frac{d}{dt}
\big\| e^{\lr{D}_{\ell}^{\rho}(T-a t)}&u\big\|^2=-2a\|\lr{D}_{\ell}^{\rho/2}v\|^2
\\
+\
((i{\tilde A} &+
{\tilde B})\,v,v)\ +\
(v,(i{\tilde A}+
{\tilde B})\,v)
\ +\
({\tilde f}, v)\ +\
(v,{\tilde f}).
\end{split}
\end{equation}

Since $i{\tilde A}+{\tilde B}\in {\tilde S}^1$ one has $
 |((i{\tilde A}+{\tilde B})\,v,v)+(v,(i{\tilde A}+{\tilde B})\,v)|\leq C\|\lr{D}_{\ell}^{1/2}v\|^2$ 
so
 \[
 \|v(t)\|^2
 \ \leq\
  \|v(0)\|^2\ +\
  C\int_0^t\|\lr{D}_{\ell}^{1/2}v\|^2ds
  \ +\
  2\int_0^t\|v\|\|{\tilde f}\|^2d\tau.
 \]
Replacing $v$ by $\lr{D}_{{\ell}}^{(\rho-1)/2}v$ yields
 \begin{align*}
 \|\lr{D}_{\ell}^{(\rho-1)/2}v(t)\|^2
 &\leq \|\lr{D}_{\ell}^{(\rho-1)/2}v(0)\|^2
 \\
 +\ C & \int_0^t\|\lr{D}_{\ell}^{\rho/2}v\|^2d\tau
 \ +\ 
 2\int_0^t\|\lr{D}_{\ell}^{(\rho-1)/2}\|\|{\tilde f}\|^2d\tau.
 \end{align*}
On the other hand,  the reasoning leading to \eqref{eq:eneid} yields
\[
\frac{d}{dt}\big(R\, v, v\big)
\ +\
ca\|\lr{D}_{{\ell}}^{\rho/2}v\|^2
\ \leq\
 C\|\lr{D}_{{\ell}}^{(\rho-1)/2}v\|\ \|\lr{D}_{{\ell}}^{2\nu-(\rho-1)/2}{\tilde f}\|
 \,.
\]
%
%
%
If
 \begin{equation}
 \label{eq:taifu}
 (\rho-1)/2\geq -\rho+1/2+2\nu\,,
 \quad
{\rm equivalently},
\quad
 \rho\geq (2+4\theta)/(3+4\theta)
 \end{equation}
then we can control $(Rv,v)$ taking \eqref{eq:Rsita} into account.  Since $ (2+4\theta)/(3+4\theta)\geq (1+3\theta)/(2+3\theta)$ and $2\nu-(\rho-1)/2\geq 0$ if \eqref{eq:taifu} is verified then 
\[
\|\lr{D}_{{\ell}}^{(\rho-1)/2}v(t)\|^2
\ \leq\
 c\,\|\lr{D}_{{\ell}}^{\nu}v(0)\|^2
 \ +\
 2\,CM\int_0^t\|\lr{D}_{{\ell}}^{2\nu-(\rho-1)/2}{\tilde f}||   \ d   \tau
\]
 where $M:=\sup_{0\leq \tau\leq t}\|\lr{D}_{{\ell}}^{(\rho-1)/2}v(\tau)\|$. Repeating the same arguments proving Proposition \ref{pro:itoyu:a} yields the following
  alternative {\it a priori} estimate.

 \begin{proposition}
 \label{pro:itoyu:b}
 If $\rho\geq (2+4\theta)/(3+4\theta)$,
 then there exist $T>0$, $a>0$, and  ${\ell}_0>0$, so that for  any $T_1<T$ there is $C>0$ such that  for all $u$ so that
 $e^{\lr{D}^{\rho}_{\ell}(T-at)}
\partial_{t ,x}^\gamma
u\in L^1([0,T]\,;\,H^{2\nu-(\rho-1)/2}(\RR^d))$ for $|\gamma|\le 1$,
\begin{align*}
\|\lr{D}_{{\ell}}^{(\rho-1)/2}e^{\lr{D}_{\ell}^{\rho}(T-a t)}u\|
\leq C\|  &  \lr{D}_{{\ell}}^{\nu}e^{T\lr{D}_{\ell}^{\rho}}u(0)\|
\\
&+C\int_0^t\big\|\lr{D}_{{\ell}}^{2\nu-(\rho-1)/2}e^{\lr{D}_{\ell}^{\rho}(T-a \tau)}Lu\big\|\ d\tau
\end{align*}
for  $0\leq t\leq T_1/a$ and ${\ell}\geq {\ell}_0$, where $\nu=\theta(\rho-1)$.
\end{proposition}
%
 

{\bf Step IV.  Uniform estimates for regularized equations.}  
Take $\chi(x)\in C_0^{\infty}({\mathbb R}^d)$ 
that  is 
equal to $1$ on a  neighborhood of  $x=0$ and such that $|\chi(x)|\leq 1$.
Consider regularized operator
\[
L^{h}\ :=\
\partial_t-\chi(h D)\big(iA(t,x,D)
\ +\
B(t,x)\big)\chi(h D)
\ :=\ 
\partial_t-iA^{h}-B^{h}.
\]
Denote 
\[
{\tilde L}^h
\ :=\
e^{\lr{D}_{\ell}^{\rho}(T-at)} L^he^{-\lr{D}_{\ell}^{\rho}(T-at)}
\ :=\
\partial_t-i{\tilde A}^h-{\tilde B}^h
\,,
\quad
{\rm and}
\quad
{\tilde L}^0:={\tilde L}\,.
\]
Denote $\chi_{h}(D):=\chi(h D)$  so
\begin{align*}
{\tilde A}^{h}
\ &=\ 
e^{\lr{D}_{\ell}^{\rho}(T-at)}\chi_{h} A \chi_h e^{-\lr{D}_{\ell}^{\rho}(T-at)}=\chi_h {\tilde A} \chi_h\, ,
\\
{\tilde B}^{h}\ &=\
e^{\lr{D}_{\ell}^{\rho}(T-at)}\chi_{h}B\chi_{h}e^{-\lr{D}_{\ell}^{\rho}(T-at)}
\ =\
\chi_{h}{\tilde B}\chi_{h}\, .
\end{align*}
Note that $\chi_{h}(\xi)\in {\tilde S}^0$ uniformly in $0<h\leq \ell^{-1}$ because $\lr{\xi}_{\ell}\leq C|\xi|$
on the support of $\nabla_{\xi}\chi(h \xi)$.

Recall that ${\tilde A}=H_N+K$ with $K$ in \eqref{eq:torii}. Since $H_N\in {\tilde S}^1$  it follows that
\[
\chi(h\xi)\, \#\, H_N \, \#\ \, \chi(h \xi)
\ =\
\chi^2(h\xi)H_N\ +\
K_1^h
\]
where $K_1^h\in {\tilde S}^0$, uniformly in $0<h\leq\ell^{-1}$. 
 It is clear that $\chi_h\#{\tilde B}\#\chi_h\in {\tilde S}^0$ and 
 $\chi_h\#K\#\chi_h\in {\tilde S}^{\max\{\rho-N(1-\rho)\,,\,-1+\rho\}}$ 
 uniformly in $0<h \leq \ell^{-1}$. 

From here  on the  pseudodifferential calculus is understood to be uniform in $0<h\leq \ell^{-1}$. Denote 
$H^{h}_N=\chi^2(h \xi)H_N$ so that
\[
H^{h}_N(\ell,\tau,t,x,\xi)=\chi_{h}^2(h \xi)\ 
\lr{\xi}_\ell\ {\mathcal H}_N\Big(t\,,\, x\,,\, \xi/\lr{\xi}_\ell\,;\, \tau\rho\lr{\xi}_\ell^{\rho-1} \Big).
\]
Choosing $s\chi_h^2\lr{\xi}_\ell$, $\tau\rho\lr{\xi}_\ell^{\rho-1}$ ($\tau>0$), $\xi/\lr{\xi}_\ell$ for $s$, $\epsilon$, $\xi$ in \eqref{eq:Matexp} 
yields
\[
\frac{ \tau^{\theta}  }
{C\, \lr{\xi}_\ell^{\theta(1-\rho)} \, e^{cs\tau \chi_h^2\lr{\xi}_\ell^{\rho}}}
\ \leq\ 
\big\| e^{isH_N^h(\ell,\tau,t,x,\xi)} \big\|
\ \leq\
\frac{C\, \lr{\xi}_\ell^{\theta(1-\rho)}   \, e^{cs\tau\chi_h^2\lr{\xi}_\ell^{\rho}}  }
{\tau^{\theta}   }\ .
\]
Define
\[
M^h\ :=iH_N^h(\ell,\tau,t,x,\xi)-a\lr{\xi}^{\rho}_{\ell},\qquad \tau=T-at
\]
and the corresponding symmetrizer 
\[
R^h(t,x,\xi)\ :=\ a\,
\int_0^{\infty}\lr{\xi}_{\ell}^{\rho}\ \big(e^{sM^h(\ell,\tau,t,x,\xi)}\big)^*\
\big(e^{sM^h(\ell,\tau,t,x,\xi)}\big)\ ds\ .
\]
Since  $\|e^{sM^h}\|=e^{-as\lr{\xi}^{\rho}_{\ell}}\|e^{sM^h}\|$ and $0\leq \chi^2_h\leq 1$ one  has
\[
\tau^{\theta}\lr{\xi}_\ell^{-\theta(1-\rho)}e^{-c_1as\lr{\xi}_\ell^{\rho}}/C
\ \leq\
 \big\|e^{sM^h}\big\|
 \ \leq\
  C \tau^{-\theta}\lr{\xi}_\ell^{\theta(1-\rho)}\ e^{-c_2a s\lr{\xi}_\ell^{\rho}}
\]
with  $c_i>0$, $C>0$ independent of $0<h\leq \ell^{-1}$, $\ell$ and $a$.  Since
\[
|\partial_{\xi}^{\alpha}\partial_x^{\beta}M^h|\leq C_{\alpha\beta}\lr{\xi}_{\ell}^{1-|\alpha|}
\]
uniformly in $0<h\leq \ell^{-1}$ the estimates for $R^h$ are
exactly the same as those for $R$,  so one has \eqref{eq:estRR} with $C_{\alpha\beta}$ independent of $0<h\leq \ell^{-1}$, $\ell$, $x$, $\xi$ and $a$. Repeating the same arguments proving Proposition \ref{pro:itoyu:a} 
proves the following.

{\sl
If  $\rho\geq (1+6\theta)/(2+6\theta)$ then there exist $T>0$, $a>0$ and ${\ell}_0>0$ such that   for any $T_1<T$ one can find $C>0$ such that for all $v$ so that $v\in C^1([0,T]; H^{3\nu}(\RR^d))$
\begin{equation}
\label{eq:existence:a1}
\|\lr{D}_{{\ell}}^{-\nu}(t)v\|
\ \leq\
 C\|\lr{D}_{{\ell}}^{\nu}v(0)\|
\ +\
C\int_0^t\|\lr{D}_{{\ell}}^{3\nu}{\tilde L}^hv(\tau)\|\,d\tau
\end{equation}
for  $0\leq t\leq T_1/a$ and ${\ell}\geq {\ell}_0$  where $C$ is independent of $\ell$ and $0<h\leq \ell^{-1}$.
}

\vskip.1cm
{\bf Step V.  Construction of solution.}   Next solve
\begin{equation}
\label{eq:hODE}
{\tilde L}^hv^h\ =\ 
(\partial _t-i{\tilde A}^h-{\tilde B}^h)v^h\ =\
{\tilde f},\qquad v^h(0)\ =\ {\tilde g}.
\end{equation}
Since $i{\tilde A}+{\tilde B}\in C({\mathbb R};{\tilde S}^1)$ and $\chi_h\in S^{-1}$
with $h$-dependent bound, it follows that 
 $i{\tilde A}^h+{\tilde B}^h\in C({\mathbb R};{\tilde S}^0)$ so is
 bounded  from $H^k({\mathbb R}^d)$ to $H^k({\mathbb R}^d)$ for any $k\in {\mathbb R}$. Therefore 
for any ${\tilde g}\in H^k(\RR^d)$ and ${\tilde f}\in L^1_{loc}({\mathbb R};H^k(\RR^d))$ there exists a unique
solution $v^h\in C^1({\mathbb R};H^k(\RR^d))$ to the linear ordinary differential equation
\eqref{eq:hODE}. 

Assume 
\[
{\tilde f}=e^{\lr{D}_{\ell}^{\rho}(T-at)}f\in L^1([0,T'];H^{3\nu}(\RR^d)),\quad {\tilde g}=e^{T\lr{D}_{\ell}^{\rho}}g\in H^{3\nu}(\RR^d)\,.
\]
Denote $T^\prime:=T_1/a$ and the corresponding solutions to 
\eqref{eq:hODE} by $v^h\in C^1([0,T'];H^{3\nu}(\RR^d))$. Then \eqref{eq:existence:a1} yields
\[
\|\lr{D}_{\ell}^{-\nu}v^h(t)\|
\ \leq\
 C\,\|\lr{D}_{\ell}^{\nu}{\tilde g}\|
 \ +\
 C\int_0^t\|\lr{D}_{\ell}^{3\nu}{\tilde f}\|\ d\tau
\]
for $0\leq t \leq T'$.  Therefore $\{v^h\}$ is bounded in $L^{\infty}([0,T'];H^{-\nu})$. Since $L^{\infty}([0,T];H^{-\nu}(\RR^d))$ is the dual of $L^1([0,T'];H^{\nu}(\RR^d))$,
  one can choose a subsequence
  (still denoted by $\{v^h\}$)  weak$^*$ convergent  in $L^{\infty}([0,T'];H^{-\nu}(\RR^d))$ to $v$.
      It is easy to see that $\chi(h D)v^h$ converges to $v$ weakly in $L^{\infty}([0,T'],H^{-\nu}(\RR^d))$. Since $i{\tilde A}+{\tilde B}$ maps $L^{\infty}([0,T'],H^{-\nu}(\RR^d))$ into $L^{\infty}([0,T'];H^{-\nu-1})$ then $\chi(h D)(i{\tilde A}+{\tilde B})\chi(h D)v^h$ converges to $(i{\tilde A}+{\tilde B})v$  weak$^*$ in $L^{\infty}([0,T'];H^{-\nu-1}(\RR^d))$. Since it is clear that 
\[
\int_0^{T'}(\partial_tv^h,\phi)dt\to -\int_0^{T'}(v,\partial_t\phi)dt
\]
for any $\phi\in C_0^{\infty}((0,T')\times {\mathbb R}^d)$ 
it follows that $v$ satisfies ${\tilde L}v={\tilde f}$ and $v(0)={\tilde g}$, that is 
\[
e^{\lr{D}_{\ell}^{\rho}(T-at)}Le^{-\lr{D}_{\ell}^{\rho}(T-at)}v\ =\
{\tilde f}\ =\
e^{\lr{D}_{\ell}^{\rho}(T-at)}f\,.
\]
The equation ${\tilde L}v={\tilde f}$ yields $\partial_tv\in L^{\infty}([0,T'];H^{-\nu-1}(\RR^d))$ which implies $v\in C([0,T'];H^{-\nu-1}(\RR^d))$. 
Since $u=e^{-\lr{D}_{\ell}^{\rho}(T-at)}v$ we conclude that
\[
Lu=f,\quad u(0)=g,\quad (t,x)\in (0,T')\times {\mathbb R}^d\, .
\]
This completes the proof of existence of a solution $u$ with
$e^{\lr{D}_{\ell}^{\rho}(T-at)}u=v\in L^{\infty}([0,T'],H^{-\nu}(\RR^d))$.

\vskip.1cm
{\bf Step VI.  Proof of uniqueness.}  Suppose that $u$
is a solution with vanishing data $f,g$.
Define $0\le t_1  \le T_0$ so that 
$u$ vanishes on $[0,t_1]\times\RR^d$ but does not
vanish on $[0,t_1+\eps)\times\RR^d$ for any $\eps>0$.
Need to show that $t_1=T_0$.
Suppose that $t_1<T_0$.

Using Remark \ref{rem:Thm1.1} applied to the adjoint operator
with time reversed,  
choose $0< \ut\le T_0-t_1$ and $C\gg 1$ so that 
for $F(t,x)$ compactly supported in $x$ and satisfying
\begin{equation}
\label{eq:charlie}
\int_0^{T_0} \Big( \int_{\RR^d}
\big| \hat F(t,\xi)\big|^2\ 
e^{C\langle \xi\rangle^{1/s} }\ d\xi
\Big)^{1/2}
\ dt\ \ < \ \infty
\end{equation}
the adjoint problem 
$$
L^* w \ =\ F \quad {\rm on}\quad (t_1,t_1+\ut)\times\RR^d\,,
\qquad
w\big|_{t=t_1+\ut}\ =\ 0
$$
has a solution in $C([t_1,t_1+\ut]\,;{\mathcal G}_0^{s}(\RR^d))$.  

Both $u$ and $w$ being compactly supported in $x$ belong to 
$H^1((t_1,t_1+\ut)\times\RR^d)$ so integration by parts shows that
with integrals over $(t_1,t_1 + \ut)\times\RR^d$,
$$
\int\int
(u,F)\ dxdt
\ =\ 
\int\int
(u,L^* w)\ dxdt
\ =\ 
\int\int
(Lu,w)\ dxdt
\ =\ 0\,,
$$
where the initial conditions $u(t_1)=w(t_1+\ut)=0$ eliminate the boundary
contributions from $t=t_1,t_1+\ut$.

Since the set of such $F$ satisfying
\eqref{eq:charlie} is dense in $L^2([t_1,t_1+\ut]\times\RR^d)$ it follows
that $u=0$ on $[t_1,t_1+\ut]\times\RR^d$.
Therefore $u$ vanishes on $[0,t_1+\ut]\times\RR^d$  violating the choice of 
$t_1$.  Thus one must have $t_1=T_0$ proving uniqueness.

\vskip.1cm
{\bf Step VII.  Proof of continuity in time.}  
Compute $\partial_tu=e^{-\lr{D}^{\rho}_{\ell}(T-at)}(a\lr{D}_{\ell}^{\rho}v+\partial_tv)$.
Since
 $a\lr{D}^{\rho}_{\ell}v+\partial_tv\in L^{\infty}([0,T'];H^{-\nu-1}(\RR^d))$ it follows that  for any $0<c<T-T_1$ 
\[
\int \Big(|{\hat u}(t,\xi)|^2
\ +\ |\partial_t{\hat u}(t,\xi)|^2\Big)\ e^{2c\lr{\xi}^{1/s}}\ d\xi
\ \in \ L^\infty([0,T^\prime])\,.
\]
This implies that $u$ is continuous with values in ${\mathcal G}^s_0(\RR^d)$.

This completes the proof of Theorem \ref{thm:bronshtein1}.
\hfill
\qed

\section{Theorem \ref{thm:lowreg}, examples and proof }

Before the details of the proof of Theorem \ref{thm:lowreg}, 
we present two examples that illustrate the 
conclusion.

\begin{example}\rm 
If  $A(t,x,\xi)$ is uniformly diagonalizable we can choose $\theta=0$ and Theorem \ref{thm:lowreg} holds with $1<s\leq 2/(2-\kappa)$.   This is weaker than the 
sharp condition $s< 1/(1-\kappa)$ of \cite{CDS}  for $u_{tt}-a(t)u_{xx}=0$, $a> 0$. 
\end{example}
\begin{example}
\rm
If \eqref{eq:Matexp} holds with $\theta=\mu-1$, $2\leq \mu\leq m$ then 
the constraints on $s$
read
\[
1<s\leq  \min{\Big\{(3\mu-1)/(3\mu-1-\kappa)\ ,\ (4\mu-1)/(4\mu-2)\Big\}}\,.
\]
This is far from the sharp bound $s<1+\kappa/2$ of \cite{CJS} in the 
case $u_{tt}-a(t) u_{xx}=0$ with $a\ge 0$ and $\theta=m-1=1.$ 
\end{example}

{\bf Proof of Theorem \ref{thm:lowreg}.}  We present only the {\it a priori}
estimate.  Existence and uniqueness then follow as in the preceding
section.  We follow the argument in \cite{Ni:1} (also \cite{Jann}).
 By hypothesis,
\[
|\partial_x^{\alpha}(A_j(t,x)-A_j(\tau,x))|\ \leq\  C A^{|\alpha|}|\alpha|!^s|t-\tau|^{\kappa}\,,
\qquad
0<\kappa\leq 1\,.
\]
Choose $\chi(s)\in C_0^{\infty}(\RR)$ such that $\chi(s)=\chi(-s)$ with
$
\int \chi(s)ds=1$.
Define, with $0<\delta$ to be chosen later,
\[
{\tilde R}(t,x,\xi)\ :=\
\lr{\xi}_{{\ell}}^{\delta}\int R(\tau,x,\xi)  \ \chi\big((t-\tau) \,  \lr{\xi} ^{\delta}_{{\ell}}  \big)\ d\tau.
\] 
Since 
$|\partial_{\xi}^{\alpha}\chi((t-\tau)\lr{\xi}_{{\ell}}^{\delta})|\leq C_{\alpha}\lr{\xi}_{{\ell}}^{-|\alpha|}$ 
Theorem \ref{thm:symbol} implies  that
${\tilde R}\in S\big(\lr{\xi}_{{\ell}}^{2\nu}, G\big)$ 
with $G$ from \eqref{eq:metricG} . It is clear that ${\tilde R}\geq c\,\lr{\xi}_{{\ell}}^{-2\nu}$. 
\begin{lemma}
\label{lem:RH} $R(t)-R(\tau)\in S\big(|t-\tau|^{\kappa}\lr{\xi}_{{\ell}}^{3\nu+1-\rho},G\big)$ uniformly in $t$, $\tau$.
That is, for all $\alpha$, $\beta$,
\begin{align*}
|\partial_x^{\beta}\partial_{\xi}^{\alpha}(R(t)-R(\tau))|
\leq C_{\alpha\beta}a^{-|\alpha+\beta|}|t-\tau|^{\kappa}\lr{\xi}_{\ell}^{3\nu+1-\rho}\lr{\xi}_{\ell}^{(1-\rho+\nu)|\beta|-(\rho-\nu)|\alpha|}.
\end{align*}
\end{lemma}
{\bf Sketch of proof of Lemma.} It suffices to repeat  arguments similar to those proving Theorem \ref{thm:symbol}. 
\hfill
\qed
\vskip.2cm

Since
\[
{\tilde R}(t)-R(t)
\ =\
\lr{\xi}_{{\ell}}^{\delta}\int 
\big(R(\tau)-R(t)\big)
\ 
\chi((t-\tau)
\ \lr{\xi}^{\delta}_{{\ell}})\ d\tau,
\]
Lemma \ref{lem:RH} implies that 
\[
{\tilde R}(t)-R(t) \ \in\
 S \big( \lr{\xi}_{{\ell}}^{3\nu +1-\rho-\kappa\delta} \, ,G \big).
\]
Similarly,
\begin{align*}
\partial_t{\tilde R}(t)&=\lr{\xi}_{{\ell}}^{2\delta}\int R(\tau,x,\xi)
\ \chi' \big( (t-\tau)\lr{\xi} ^{\delta}_{{\ell}}  \big)\ d\tau
\\
&=\lr{\xi}_{{\ell}}^{2\delta}\int\big(R(\tau)-R(t)\big)
\ \chi' \big((t-\tau)\lr{\xi}_{{\ell}}^{\delta}  \big)\ d\tau
\end{align*}
implies that
\begin{equation}
\label{eq:dtR}
\partial_t{\tilde R}(t)\in S \big( \lr{\xi}_{{\ell}}^{3\nu +1-\rho+\delta-\kappa\delta},G \big).
\end{equation}

With 
 ${\tilde K}=i{\tilde A}+{\tilde B}$ and ${\tilde  f}=e^{\lr{D}_{\ell}^{\rho}(T-a t)}f$, one has
\begin{align*}
\frac{d}{dt}({\tilde R}\ e^{\lr{D}_{\ell}^{\rho}(T-a t)}u,e^{\lr{D}_{\ell}^{\rho}(T-a t)}u)
=
(\partial_t{\tilde R}\ v,v) 
+({\tilde R}({\tilde K}-a\lr{D}_{\ell}^{\rho})v,v)\\
+({\tilde R}v,({\tilde K}-a\lr{D}_{\ell}^{\rho})v)
+({\tilde R}{\tilde f}, v)+({\tilde R}v,{\tilde f}).
\end{align*}
Adding and subtracting two terms, the right hand side is equal to
\begin{align*}
&(\partial_t{\tilde R}v,v)+(R({\tilde K}-a\lr{D}_{\ell}^{\rho})v,v)+(Rv,({\tilde K}-a\lr{D}_{\ell}^{\rho})v)+
\\
&(({\tilde R}-R)({\tilde K}-a\lr{D}_{\ell}^{\rho})v,v)
+(({\tilde R}-R)v,({\tilde K}-a\lr{D}_{\ell}^{\rho})v)
+({\tilde R}{\tilde f}, v)+({\tilde R}v,{\tilde f}).
\end{align*}
For the terms 
\[
({\tilde R}{\tilde f}, v),\;\;(R({\tilde K}-a\lr{D}_{\ell}^{\rho})v,v)+(Rv,({\tilde K}-a\lr{D}_{\ell}^{\rho})v),\;\;
({\tilde R}{\tilde f}, v)+({\tilde R}v,{\tilde f})
\]
 use the same estimates as in Section \ref{Sec:shomei}.
For the other terms use \eqref{eq:dtR}, ${\tilde R}-R\in {\tilde S}^{3\nu+1-\rho-\kappa\delta}_{\rho-\nu,1-\rho+\nu}$ and ${\tilde K}-a\lr{\xi}_{{\ell}}^{\rho}\in {\tilde S}^1$ to find the pair of estimates,
\begin{align*}
\big|(({\tilde R}-R)({\tilde K}-a\lr{D}_{\ell}^{\rho})v,v)\big|+\big|(({\tilde R}-R)v, & ({\tilde K}-a\lr{D}_{\ell}^{\rho})v) \big|
\\
&\leq C\big\|\lr{D}_{{\ell}}^{(3\nu+2-\rho-\kappa\delta)/2}v\big\|^2\,,
\end{align*}
\[
|(\partial_t{\tilde R}v,v)|\leq C\|\lr{D}_{{\ell}}^{(3\nu+1-\rho+(1-\kappa)\delta)/2}v\|^2.
\]
If $3\nu+2-\rho-\kappa\delta\leq \rho$ and $3\nu+1-\rho+(1-\kappa)\delta\leq \rho$ then
both terms are bounded by $\|\lr{D}_{{\ell}}^{\rho/2}v\|^2$ and can  be absorbed in a Gronwall
estimate.  With $\kappa$ and $\nu$ fixed the region in the 
$\delta,\rho$ plane described by the two constraints
is bounded below by a pair of lines as the figure.
\begin{center}
\includegraphics[height=2.2cm, keepaspectratio]{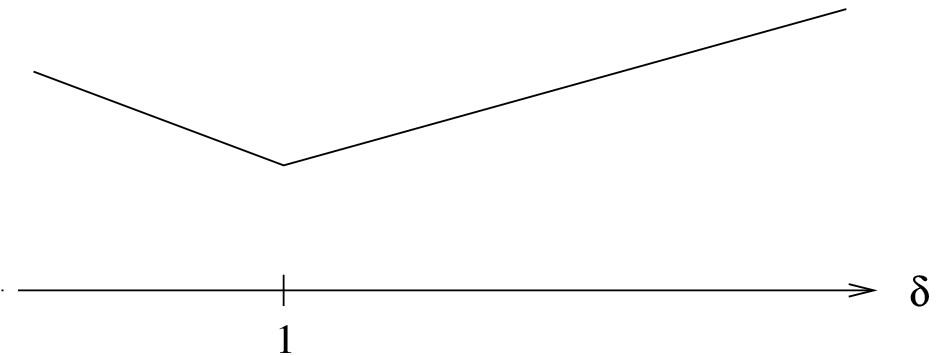}
\end{center}
The minimal value of $\rho$ 
satisfying
 the constraints occurs at $\delta=1$ independent
of $\kappa$ and $\nu$ and yields
\[
\rho\ \geq
\
\frac{ 3\theta+2-\kappa}
{3\theta+2}\,.
\]
The desired {\it a priori} estimate follows.
\hfill\qed

\section{Appendix.  The conjugation Proposition \ref{pro:weyl:1}}
\label{Sec:composition}

\begin{lemma}
\label{lem:b_1}
Let $a(x,\xi)\in {\tilde S}_{(s)}^m$ and assume $\partial_x^{\alpha}a(x,\xi)=0$ outside $|x|\leq R$ with some $R>0$ if $|\alpha|>0$. Set
\[
e^{\tau\lr{D}_{{\ell}}^{\rho}}\ a(x,D)\ e^{-\tau\lr{D}_{{\ell}}^{\rho}}\ =\ b(x,D)
\]
where $\tau\in\RR$ then $b(x,\xi)$ is given by
\begin{equation}
\label{eqn:weyl:7}
b(x,\xi)\ =\ 
\int e^{-iy\eta}\
e^{\tau\lr{\xi+\frac{\eta}{2}}^{\rho}_{\ell}-\tau\lr{\xi-\frac{\eta}{2}}^{\rho}_{\ell}}\
a(x+y,\xi)\
dyd\eta.
\end{equation}
\end{lemma}
{\bf Proof:} Write $\phi(\xi)=\tau\lr{\xi}_{{\ell}}^{\rho}$ and insert  $v=e^{-\phi(D)}u(y)=\int e^{iy\zeta-\phi(\zeta)}{\hat u}(\zeta)d\zeta$ 
into
\[
e^{\phi(D)}\ a(x,D)\ v\ =
\int e^{i(x\xi-z\xi+(z-y)\eta)}\
e^{\phi(\xi)}\
 a\Big(\frac{z+y}{2},\eta\Big)\
 v(y)\
 dyd\eta dzd\xi
\]
to get
\[
e^{\phi(D)}\ a(x,D)\ e^{-\phi(D)}u
=\int e^{ix\zeta}
\ I(x,\zeta,\mu)
\
{\hat u}(\zeta)\
d\zeta
\]
where
\[
I\ =\
\int e^{i(x\xi-z\xi+(z-y)\eta+y\zeta-x\zeta)}\
e^{\phi(\xi)}\
a\Big(\frac{z+y}{2},\eta\Big)\
e^{-\phi(\zeta)}\
dyd\eta dz d\xi.
\]
The change of variables ${\tilde z}=(y+z)/2$, ${\tilde y}=(y-z)/2$ yields
\begin{align*}
I
\ &=\
2^n\int e^{i{\tilde y} (\xi-2\eta+\zeta)}\ d{\tilde y}\
 \int e^{-i({\tilde z}-x) (\xi-\zeta)}\
 e^{\phi(\xi)}\
 a({\tilde z},\eta)\
 e^{-\phi(\zeta)}\
 d\eta d{\tilde z} d\xi
 \\
\ &=\
2^n\int e^{-2i({\tilde z}-x)(\eta-\zeta)}e^{\phi(2\eta-\zeta)}\
a({\tilde z},\eta,\mu)\
e^{-\phi(\zeta)}\
d\eta d{\tilde z}\\
\ &=\
\int e^{-i{\tilde z} \eta}\
e^{\phi(\sqrt{2}\eta+\zeta)-\phi(\zeta)}\
a\Big(x+\frac{{\tilde z}}{\sqrt{2}},\zeta+\frac{\eta}{\sqrt{2}}\Big)\
d\eta d{\tilde z}
\end{align*}
and then
\[
e^{\phi(D)}\ a(x,D)\ e^{-\phi(D)}\ u
\ =\
\int e^{i(x-y)\xi}\
p(x,\xi)\
u(y)\
dyd\xi
\ =\ 
{\rm Op}^0(p)u
\]
with
\begin{equation}
\label{eq:atu:b}
p(x,\xi)\ =\
\int e^{-iy\eta}\
e^{\phi(\xi+\sqrt{2}\eta)-\phi(\xi)}\
a\Big(x+\frac{y}{\sqrt{2}},\xi+\frac{\eta}{\sqrt{2}}\Big)
\ dyd\eta.
\end{equation}
Here we remark ${\rm Op}^0(p)=b(x,D)$ with $b(x,\xi)$ given by
\begin{equation}
\label{eq:atu}
b(x,\xi)\ =\
\int e^{iz\zeta}\
 p\Big(x+\frac{z}{\sqrt{2}},\xi+\frac{\zeta}{\sqrt{2}}\Big)\
 dzd\zeta.
\end{equation}
Indeed we see
\begin{align*}
b(x,D)u\ &=\
\int e^{i(x-y)\xi}\
b(\frac{x+y}{2},\xi)\
u(y)\
dyd\xi
\\
\ &=\
\int e^{i(x\xi-y\xi+z\zeta)}\
p\Big(\frac{x+y}{2}+\frac{z}{\sqrt{2}},\xi+\frac{\zeta}{\sqrt{2}}\Big)\
u(y)\
dyd\xi dzd\zeta\\
\ &=\
\int e^{i((x-y-z)\xi+z\zeta)}\
p\Big(\frac{x+y+z}{2},\zeta\Big)\
u(y)\
dyd\xi dz d\zeta\\
\ &=\
\int e^{iz\zeta}\
p(x,\zeta)\
u(x-z)\
dzd\zeta\ =\ {\rm Op}^0(p)u.
\end{align*}
Insert \eqref{eq:atu:b} into \eqref{eq:atu} to get
\begin{align*}
b(x,\xi)\ =\ 
\int & e^{i(z\zeta-y\eta)}\
e^{\phi(\sqrt{2}\eta+\xi+\frac{\zeta}{\sqrt{2}})-\phi(\xi+\frac{\zeta}{\sqrt{2}})}\
a\ dyd\eta dzd\zeta,
\end{align*}
$$
a \ =\ a\Big(x+\frac{z+y}{\sqrt{2}},\xi+\frac{\eta+\zeta}{\sqrt{2}}\Big).
$$
The change of variables
\[
{\tilde z}=\frac{z+y}{\sqrt{2}},\quad
{\tilde y}=\frac{y-z}{\sqrt{2}},\quad
{\tilde \zeta}=\frac{\zeta+\eta}{\sqrt{2}},\quad
{\tilde \eta}=\frac{\eta-\zeta}{\sqrt{2}}
\]
gives
\begin{align*}
b(x,\xi)\ &=\
\int e^{-i({\tilde z} {\tilde \eta}+{\tilde y}{\tilde \zeta})}\
e^{\phi(\frac{3{\tilde \zeta}}{2}+\xi+\frac{{\tilde \eta}}{2})-\phi(\xi+\frac{{\tilde \zeta}}{2}-\frac{{\tilde \eta}}{2})}\
a(x+{\tilde z},\xi+{\tilde \zeta})\
d{\tilde y}  d{\tilde \eta} d{\tilde z} d{\tilde \zeta}
\\
\ &=\
\int e^{-i{\tilde z}{\tilde\eta}}\
e^{\phi(\xi+\frac{{\tilde \eta}}{2})-\phi(\xi-\frac{{\tilde \eta}}{2})}\
a(x+{\tilde z},\xi)\
d{\tilde z}d{\tilde \eta}
\end{align*}
proving \eqref{eqn:weyl:7}.
\hfill\qed
\vskip.2cm
{\bf Proof of Proposition \ref{pro:weyl:1}}. Insert
\begin{align*}
a(x+y,\xi)=\sum_{|\alpha|\leq  k} & \frac{1}{\alpha !}D_x^{\alpha}a(x,\xi)(iy)^{\alpha}
\\
&+\sum_{|\alpha|=k+1}\frac{k+1}{\alpha !}(iy)^{\alpha}\int_0^1(1-\theta)^{k}D_x^{\alpha}a(x+\theta y,\xi)d\theta
\end{align*}
into \eqref{eqn:weyl:7} to get
\begin{equation}
\label{eq:okazaki}
\begin{split}
&b(x,\xi)=\sum_{|\alpha|\leq k}\frac{1}{\alpha !}\int e^{-iy\eta}e^{\phi(\xi+\frac{\eta}{2})-\phi(\xi-\frac{\eta}{2})}D_x^{\alpha}a(x,\xi)(iy)^{\alpha}dyd\eta\ +\
\\
\sum_{|\alpha|=k+1}&\frac{k+1}{\alpha !}\int e^{-iy\eta}e^{\phi(\xi+\frac{\eta}{2})-\phi(\xi-\frac{\eta}{2})}(iy)^{\alpha}dyd\eta
 \int_0^1(1-\theta)^{k}D_x^{\alpha}a(x+\theta y,\xi)d\theta.
\end{split}
\end{equation}
Since $e^{-iy\eta}(iy)^{\alpha}=(-\partial_{\eta})^{\alpha}e^{-iy\eta}$ the first term on the right-hand side of \eqref{eq:okazaki} is
\begin{equation}
\label{eqn:weyl:8}
\sum_{|\alpha|\leq  k}\frac{1}{\alpha !}\partial_{\eta}^{\alpha}e^{\phi(\xi+\frac{\eta}{2})-\phi(\xi-\frac{\eta}{2})}\Big|_{\eta=0}D_x^{\alpha}a(x,\xi).
\end{equation}
Note that $\partial_{\eta}^{\alpha}e^{\phi(\xi+\frac{\eta}{2})-\phi(\xi-\frac{\eta}{2})}\Big|_{\eta=0}$ is a linear combination of
\[
\partial_{\xi}^{\alpha_1}\phi(\xi)\cdots \partial_{\xi}^{\alpha_s}\phi(\xi),\qquad \sum_{j=1}^s \alpha_j=\alpha,
\quad|\alpha_j|\geq 1.
\]
Divide the linear combination into two parts; the sum over $\sum \alpha_j=1$, $|\alpha_j|=1$ and the remaining sum called $r$. If $|\alpha_j|\geq 2$ for some $j$ then $s\leq |\alpha|-1$ and hence $s\rho-|\alpha|\leq -(1-\rho)|\alpha|-\rho\leq -2+\rho$ so that
\[
\partial_{\xi}^{\alpha_1}\phi(\xi)\cdots \partial_{\xi}^{\alpha_s}\phi(\xi)\in {\tilde S}^{-2+\rho }.
\]
Then \eqref{eqn:weyl:8} yields
\[
\sum_{|\alpha|\leq  k}\frac{1}{\alpha!}D_x^{\alpha}a(x,\xi)(\tau\nabla_{\xi}\lr{\xi}_{{\ell}}^{\rho})^{\alpha}\ +\
r\,,
\qquad
r\in {\tilde S}^{-1+\rho}.
\]
Define
\begin{align*}
H_{\alpha}(\xi,\eta,\mu)\ &=\
\frac{1}{\alpha!}\, \partial_{\eta}^{\alpha}e^{\phi(\xi+\frac{\eta}{2})
-\phi(\xi-\frac{\eta}{2})}\\
\ &=\
2^{-|\alpha|}\sum_{\beta+\gamma=\alpha}\frac{1}{\beta !\gamma !}\
\partial_{\xi}^{\beta}e^{\phi(\xi+\frac{\eta}{2})}\
(-\partial_{\xi})^{\gamma}e^{-\phi(\xi-\frac{\eta}{2})}
\end{align*}
where the second term on the right-hand side of \eqref{eq:okazaki} is, up to a multiplicative constant 
\begin{align*}
\sum_{|\alpha|=k+1}\int   &  e^{-iy\eta}    H_{\alpha}(\xi,\eta)dyd\eta\int_0^1(1-\theta)^{k}D_x^{\alpha}a(x+\theta y,\xi)d\theta\\
&=\sum_{|\alpha|=k+1}\int\int_0^1 e^{ix\eta}(1-\theta)^{k}H_{\alpha}(\xi,\theta\eta)d\eta d\theta\int e^{-iy\eta}D_x^{\alpha}a(y,\xi)dy.
\end{align*}
Define $E_\alpha(\eta,\xi):=\int e^{-iy\eta}D_x^{\alpha}a(y,\xi)dy$ and
\begin{equation}
\label{eqn:weyl:9}
R_k\ :=\
\sum_{|\alpha|=k+1}\int\int_0^1 e^{ix\eta}\
(1-\theta)^{k}\
H_{\alpha}(\xi,\theta\eta)\
E_\alpha(\eta,\xi)\
d\eta d\theta\,.
\end{equation}
\begin{lemma}
\label{le:weyl:5}
There is $c>0$ such that for any $\delta\in\NN^n$,
\[
|\partial_{\xi}^{\delta}E_\alpha(\eta,\xi)|
\ \leq\
 C_{\alpha\delta}\,
 \lr{\xi}_{{\ell}}^{1-|\delta|}\,
 e^{-c\lr{\eta}^{\rho}}.
\]
\end{lemma}
{\bf Proof.} Integration by parts gives
\[
{\eta}^{\nu}\partial_{\xi}^{\delta}E_\alpha(\eta,\xi)
\ =\
\int e^{-iy\eta}\ \partial_{\xi}^{\delta}D_x^{\alpha+\nu}a(y,\xi)\ dy\,.
\]
Then there exist constants $A>0$ and $C_{\delta}$ such that 
\begin{eqnarray*}
|\partial_{\xi}^{\delta}E_\alpha(\eta,\xi)|\leq C_{\delta}\lr{\xi}_{{\ell}}^{1-|\delta|} A^{|\alpha+\nu|}|\alpha+\nu|!^{s}\lr{\eta}^{-|\nu|}
\leq C_{\alpha\delta}\lr{\xi}_{{\ell}}^{1-|\delta|}A^{|\nu|}|\nu|!^s\lr{\eta}^{-|\nu|}.
\end{eqnarray*}
Choose $\nu$ minimizing $A^{|\nu|}|\nu|!^s\lr{\eta}^{-|\nu|}$, that is $|\nu|\sim e^{-1}A^{-1/s}\lr{\eta}^{1/s}$ so that $A^{|\nu|}|\nu|!^s\lr{\eta}^{-|\nu|}\lesssim e^{-s^{-1}A^{-1/s}\lr{\eta}^{1/s}}=e^{-c\lr{\eta}^{\rho}}$.
\hfill\qed
\vskip.2cm

Returning to the proof of Proposition \ref{pro:weyl:1}},
note that $H_{\alpha}(\xi,\eta)$ is a linear combination of terms
\begin{align*}
\partial_{\xi}^{\beta_1}\phi(\xi+\frac{\eta}{2})\cdots \partial_{\xi}^{\beta_s}\phi(\xi+\frac{\eta}{2})\partial_{\xi}^{\gamma_1}\phi(\xi-\frac{\eta }{2})&\cdots\partial_{\xi}^{\gamma_t}\phi(\xi-\frac{\eta}{2})
\
 e^{\phi(\xi+\frac{\eta}{2})-\phi(\xi-\frac{\eta}{2})}
 \\
 :=\ &{k}_{\beta_1,\ldots,\beta_s,\gamma_1,\ldots,\gamma_t}(\xi,\eta)
 \,e^{\phi(\xi+\frac{\eta}{2})-\phi(\xi-\frac{\eta}{2})}
\end{align*}
where $\sum \beta_j=\beta$, $\sum \gamma_j=\gamma$ and $|\beta_j|\geq 1$, $|\gamma_j|\geq 1$, $\beta+\gamma=\alpha$. Since $\lr{\xi\pm\eta/2}_{{\ell}}^r\leq C_r\lr{\xi}_{{\ell}}^r \lr{\eta}^{|r|}$ we see 
\begin{equation}
\label{eqn:weyl:10}
|\partial_{\xi}^{\delta}{k}_{\beta_1,\ldots,\beta_s,\gamma_1,\ldots,\gamma_t}(\xi,\eta)|
\ \leq\
 C_{\delta}\,\lr{\xi}_{{\ell}}^{-|\alpha|(1-\rho)-|\delta|}
 \, \lr{\eta}^{|\alpha|+|\delta|}.
\end{equation}
For some $0<\theta<1$ one has
\begin{align}
\label{eq:syuri}
\partial_{\xi}^{\alpha}\big(\phi(\xi+\frac{\eta}{2})-\phi(\xi-\frac{\eta}{2})\big)
=\sum_{j=1}^n\frac{1}{2}\eta_j\Big(\partial_{\xi}^{\alpha}\partial_{\xi_j}\phi(\xi+\frac{\theta\eta}{2})+\partial_{\xi}^{\alpha}\partial_{\xi_j}\phi(\xi-\frac{\theta\eta}{2})\Big).
\end{align}
Then $\lr{\xi\pm\theta\eta/2}_{{\ell}}^{\rho-1-|\alpha|}\leq \lr{\xi\pm\theta\eta/2}_{{\ell}}^{-|\alpha|}\leq C_{\alpha}\lr{\xi}_{{\ell}}^{-|\alpha|}\lr{\eta}^{|\alpha|}$ and 
\begin{align*}
\Big|
\partial_{\xi}^{\alpha_1}\Big[
\phi\Big(\xi+\frac{\eta}{2}\Big)-\phi\Big(\xi-\frac{\eta}{2}\Big)\Big]
&\cdots \partial_{\xi}^{\alpha_t}\Big[
\phi\Big(\xi+\frac{\eta}{2}\Big)
-\phi\Big(\xi-\frac{\eta}{2}\Big)\Big]\Big|
\\
\ &\leq\
 C_{\alpha}\,
 \lr{\xi}_{{\ell}}^{-|\alpha|}\lr{\eta}^{2|\alpha|},\quad
 \alpha=\alpha_1+\cdots+\alpha_t
\end{align*}
yield
\begin{equation}
\label{eq:velo:a}
|\partial_{\xi}^{\delta}e^{\phi(\xi+\frac{\eta}{2})
-\phi(\xi-\frac{\eta}{2})}|
\ \leq \
C_{\delta}\,
\lr{\xi}_{{\ell}}^{-|\delta|}\
\lr{\eta}^{2|\delta|}\
e^{\phi(\xi+\eta/2)-\phi(\xi-\eta/2)}
\,.
\end{equation}
Next prove that with some $c_1>0$,
\[
\big| \phi(\xi+\eta/2)-\phi(\xi-\eta/2) \big|\ \leq\  c_1|\tau|\ \lr{\eta}^{\rho}\,.
\]
 Indeed if ${\ell}+|\xi|\geq |\eta|$ then $\lr{\xi}_{{\ell}}\approx \lr{\xi\pm \theta\eta/2}_{{\ell}}$ for $|\theta|\leq 1$ hence  \eqref{eq:syuri} gives 
\begin{align*}
|\phi(\xi+\eta/2)-\phi(\xi-\eta/2)|\ &\leq\   C|\tau| \,\lr{\eta}\,\lr{\xi\pm\theta\eta/2}^{\rho-1}_{{\ell}}
\\
\ &\leq\ 
  C'\,|\tau|\,\lr{\eta}\,\lr{\xi}^{\rho-1}_{{\ell}}\ \leq \
  C''\,|\tau|\,\lr{\eta}^{\rho}.
\end{align*}
While if ${\ell}+|\xi|\leq |\eta|$ then $\lr{\xi\pm\eta/2}_{{\ell}}\leq C\lr{\eta}$ and the assertion is clear. From \eqref{eqn:weyl:10} and \eqref{eq:velo:a} we have
\begin{equation}
\label{eq:Halpha} 
|\partial_{\xi}^{\delta}H_{\alpha}(\xi,\eta)|
\ \leq\
 C_{\alpha\delta}\
 \lr{\xi}_{{\ell}}^{-|\alpha|(1-\rho)}\
 \lr{\xi}_{{\ell}}^{-|\delta|}\
 \lr{\eta}^{|\alpha|+2|\delta|}\
 e^{c_1|\tau| \lr{\eta}^{\rho}}.
\end{equation}
From Lemma \ref{le:weyl:5} and \eqref{eq:Halpha} one has
\begin{eqnarray*}
|\partial_{\xi}^{\delta}(H_{\alpha}(\xi,\eta)E_\alpha(\eta,\xi))|
\ \leq \
C_{\alpha\delta}\
 \lr{\xi}_{{\ell}}^{1-|\delta|-|\alpha|(1-\rho)}\
\lr{\eta}^{|\alpha|+2|\delta|}\
e^{-(c-c_1|\tau|)\lr{\eta}^{\rho}}
\end{eqnarray*}
where $c>0$ is the constant in Lemma \ref{le:weyl:5}. If $c-c_1|\tau|>0$ then
\[
|\partial_x^{\beta}\partial_{\xi}^{\delta}R_k(x,\xi)|
\ \leq\
 C_{\delta\beta}\,
 \lr{\xi}_{{\ell}}^{1-|\delta|-(k+1)(1-\rho)}.
\]
Since $1-(k+1)(1-\rho)=\rho-k(1-\rho)$, the assertion follows.
\hfill\qed

\end{document}